\documentclass[11pt]{article}       

\usepackage{graphicx}

\usepackage{amssymb}
\usepackage{amsthm}
\usepackage{txfonts}
\usepackage{tabularx}
\usepackage{epsfig}
\usepackage{mathdots}

\usepackage{booktabs}
\usepackage{threeparttable}

\newtheorem{proposition}{Proposition}
\newtheorem{lemma}{Lemma}
\newtheorem{theorem}{Theorem}
\newtheorem{corollary}{Corollary}
\newtheorem{conjecture}{Conjecture}

\title{Mader's Conjecture and Its Variants for Cographs}

\author{Toru Hasunuma \\ \\ 
Department of Mathematical Sciences, Tokushima University, \\
2--1 Minamijosanjima, Tokushima 770--8506 Japan}

\date{November 16, 2025}

\begin{document}

\maketitle

\begin{abstract}
The class of cographs is one of the most well-known graph classes,
which is also known to be equivalent to the class of $P_4$-free graphs.
We show that Mader's conjecture is true if we restrict ourselves to cographs, 
that is, for any tree $T$ of order $m$, every $k$-connected cograph $G$
with $\delta(G) \geq \left\lfloor \frac{3k}{2} \right\rfloor +m-1$ 
contains a subtree $T' \cong T$ such that $G-V(T')$ is still $k$-connected,
where $\delta(G)$ denotes the minimum degree of $G$.
Moreover, we show that three variants of Mader's conjecture hold for cographs, that is, 
for any tree $T$ of order $m$,
\begin{itemize}
\item every $k$-connected (respectively, $k$-edge-connected) 
cograph $G$ with $\delta(G) \geq k+m-1$ 
contains a subtree $T' \cong T$ such that $G-E(T')$ is $k$-connected 
(respectively, $k$-edge-connected),
\item every $k$-edge-connected cograph $G$ with $\delta(G) \geq k+m-[k = 1]$  
contains a subtree $T' \cong T$ such that $G-V(T')$ is $k$-edge-connected,
where we use Iverson's convention for $[k = 1]$.
\end{itemize}
We furthermore present tight lower bounds on the minimum degree of a cograph 
for the existence of disjoint connectivity keeping trees, 
a maximal connectedness keeping tree 
and a super edge-connectedness keeping tree. 

\bigskip

\end{abstract}

\section{Introduction}

Throughout this paper, a graph means a simple undirected graph.
Let $G = (V,E)$ be a graph. 
The order of $G$ is $|V(G)|$; we may simply denote it by $n_G$, i.e., $n_G = |V(G)|$. 
A graph of order 1 is called trivial.
For $v \in V(G)$, we denote by $N_G(v)$ the set of vertices adjacent to $v$ in $G$,
i.e., $N_G(v) = \{ w \in V(G)\ |\ vw \in E(G) \}$.
The degree of a vertex $v$ in $G$ is denoted by ${\rm deg}_G(v)$,
i.e., ${\rm deg}_G(v) = |N_G(v)|$.
A vertex of degree 0 (respectively, 1) is called an isolated vertex
(respectively, a leaf).
Let $\delta(G) = \min_{v \in V(G)}{\rm deg}_G(v)$ and $\Delta(G) = \max_{v \in V(G)}{\rm deg}_G(v)$.
For $S \subsetneq V(G)$, 
$G-S$ denotes the graph obtained from $G$ by deleting every vertex in $S$.  
When $S = \{s\}$, $G-S$ may be abbreviated to $G-s$.
For $F \subseteq E(G)$ and $e \in E(G)$, $G-F$ and $G-e$ are similarly defined.
A vertex-cut (respectively, edge-cut) of a connected graph $G$ is a subset $S \subsetneq V(G)$ 
(respectively, $F \subseteq E(G)$) such that $G-S$ (respectively, $G-F$) is disconnected. 
For two sets $A$ and $B$, 
$A \setminus B$ denotes the set difference $\{x\ |\ x \in A, x \not\in B \}$.
For a nonempty subset $S$ of $V(G)$, we denote by $\langle S \rangle_G$
the induced subgraph of $G$ by $S$, i.e., $\langle S \rangle_G = G - (V(G) \setminus S)$.
If for any $S \subseteq V(G)$, $\langle S \rangle_G \not\cong H$ for some graph $H$,
then $G$ is called $H$-{\it free}.
The complete graph, the path ant the cycle of order $n$ are denoted by $K_n$, $P_n$ and $C_n$,
respectively.  
For a true-or-false statement $P$, Iverson's convention $[P]$ is defined to be 1 (respectively, 0)
if $P$ is true (respectively, false).

Many graph operations have been introduced in graph theory until now.
For example, 
\begin{itemize}
\item unary operations: complement, line graph, subdivided-line graph, power, etc., 
\item binary operations: union, join, Cartesian product, lexicographic product, 
exponentiation, edge sum, etc.
\end{itemize}
These operations are actually treated in \cite{BCL,B,HF,HaS,HaE}. 
In particular, complement, union and join are the most traditional ones.
These three operations are defined as follows,
where $G_1$ and $G_2$ are assumed to be disjoint, i.e., $V(G_1) \cap V(G_2) = \emptyset$.
\begin{itemize}
\item The {\it complement} of $G$ denoted $\overline{G}$ is the graph with $V(\overline{G}) = V(G)$
and $E(\overline{G}) = \{ uv\ |\ u,v \in V(G), uv \not\in E(G) \}$.
\item The {\it union} of $G_1$ and $G_2$ denoted $G_1 \cup G_2$
is the graph with $V(G_1 \cup G_2) = V(G_1) \cup V(G_2)$
and $E(G_1 \cup G_2) = E(G_1) \cup E(G_2)$.
\item The {\it join} of $G_1$ and $G_2$ denoted $G_1 + G_2$
is the graph with $V(G_1 + G_2) = V(G_1) \cup V(G_2)$
and $E(G_1 + G_2) = E(G_1) \cup E(G_2) \cup \{ uv\ |\ u \in V(G_1), v \in V(G_2) \}$.
\end{itemize}
There are several different notations for the join in literature. 
Our notation for the join follows the texts \cite{BCL,B,HF}.
Using the union and the join, {\it cographs} are recursively defined as follows:
\begin{enumerate}
\item $K_1$ is a cograph, 
\item The union of cographs is a cograph, 
\item The join of cographs is a cograph.
\end{enumerate}
Since $G_1 + G_2 = \overline{\overline{G_1} \cup \overline{G_2}}$, 
the cographs can also be defined recursively by replacing the third condition 
with the condition that the complement of a cograph is a cograph;
in fact, the term ``cograph" is an abbreviation of the term ``complement reducible graph"
introduced in \cite{CLB}.
The class of cographs is an important subclass of the perfect graphs and
has been widely investigated so far (e.g., see \cite{ALW,BM,CPS,FNSS,J,LWB}). 
The notion of cographs itself was historically introduced by various researchers independently;
there are actually many characterizations of cographs (see \cite{CLB}).
In particular, one of the most well-known characterizations 
is that $G$ is a cograph if and only if $G$ is $P_4$-free.
From this fact, there are papers in which cographs are simply defined to be $P_4$-free graphs
without introducing the graph operations. 

The connectivity $\kappa(G)$ of a graph $G$ is the minimum number of vertices
whose removal from $G$ results in a disconnected graph or a trivial graph.
The edge-connectivity $\lambda(G)$ of $G$ is similarly defined by replacing vertex-removal
with edge-removal.
A graph $G$ is $k$-{\it connected} (respectively, $k$-{\it edge-connected}) 
if $\kappa(G) \geq k$ (respectively, $\lambda(G) \geq k$).
Note that $\kappa(K_1) = \lambda(K_1) = 0$ although $K_1$ is connected. 
Then, for convenience, the trivial graph $K_1$ is assumed to be 1-connected and 1-edge-connected 
throughout the paper.
As a fundamental property, it holds for any graph $G$ that 
$$\kappa(G) \leq \lambda(G) \leq \delta(G).$$
When $\kappa(G) = \delta(G)$ (respectively, $\lambda(G) = \delta(G))$, 
$G$ is called {\it maximally connected} (respectively, {\it maximally edge-connected}).
If every minimum edge-cut of a connected graph $G$ isolates
a vertex, i.e., for any edge-cut $F$ with $|F| = \lambda(G)$, 
there exists a vertex $v$ of degree $\delta(G)$
such that $F = \{vw \ |\ w \in N_G(v)\}$, then $G$ is called {\it super edge-connected}. 
We may define a disconnected graph $G$ to be super edge-connected if
$G$ has only one isolated vertex.

Mader posed in 2010 the following conjecture concerning the existence of a 
connectivity keeping tree.

\begin{conjecture} {\rm (Mader \cite{M1})} \label{Mader} 
For any tree $T$ of order $m$,
every $k$-connected graph $G$ with $\delta(G) \geq 
\left\lfloor \frac{3k}{2} \right\rfloor +m -1$ contains a subtree $T' \cong T$ such that
$G-V(T')$ is $k$-connected.
\end{conjecture}

Mader's conjecture generalizes a fundamental result shown in 1972 by Chartrand, Kaugers and Lick
and a well-known proposition (e.g., see \cite{BCL}) on the existence of a subtree
isomorphic to any given tree.

\begin{theorem} {\rm (Chartrand, Kaugers and Lick \cite{CKL})}
Every $k$-connected graph $G$ with $\delta(G) \geq \left\lfloor \frac{3k}{2} \right\rfloor$ 
contains a vertex $v$ such that $G-v$ is $k$-connected.
\end{theorem}

\begin{proposition} \label{tree}
For any tree $T$ of order $m$,
every graph $G$ with $\delta(G) \geq m-1$ contains a subtree $T' \cong T$.
\end{proposition}

Mader showed in \cite{M1} that the conjecture holds if $T$ is a path. 
When stating the conjecture, Mader also remarked that
Diwan and Tholiya \cite{DT} had already shown in 2009 the same statement as the conjecture for $k = 1$;
their result was actually motivated by Locke's conjecture. 
After being presented several affirmative answers to the conjecture for $k = 2$ (e.g., \cite{H0,HO,LZ}), 
the conjecture for $k = 2,3$ has been shown to be true in \cite{HL}.
Mader's conjecture remains open for $k \geq 4$. 

In this paper, we show that Mader's conjecture holds for all $k \geq 1$
if we restrict ourselves to cographs.
We in fact prove the statement with smaller lower bounds on $\delta(G)$ in several cases.
As far as we know, the class of cographs is the first graph class for which Mader's conjecture holds
for all $k \geq 1$, except for specific graphs such as complete graphs and paths. 
We also show the following two interesting statements related to Mader's conjecture. 
One is that if we replace the floor function with the ceiling function for the lower bound 
on $\delta(G)$, then we have two disjoint connectivity keeping trees when $k = \kappa(G)$. 
The other is that Proposition \ref{tree} can be extended to the existence of 
a maximal connectedness keeping tree for cographs except for $K_m$. 
In particular, in order to show the second result, we present a characterization of
the maximally connected cographs. 
\begin{itemize}
\item For any trees $T_1$ and $T_2$ of order $m$, 
every connected cograph $G$ with $\delta(G) \geq \left\lceil \frac{3\kappa(G)}{2} \right\rceil +m-1$
contains disjoint subtrees $T'_1 \cong T_1$ and $T'_2 \cong T_2$ such that 
$\kappa(G-V(T'_1) \cup V(T'_2)) = \kappa(G)$.
\item For any tree $T$ of order $m$, every maximally connected cograph $G \not\cong K_m$ 
with $\delta(G) \geq m-1$
contains $T' \cong T$ such that $G-V(T')$ is maximally connected.
\end{itemize}

As the edge-version a similar problem, 
we can consider the graph obtained from $G$ by deleting every 
edge of a subtree $T'$ instead of every vertex of $T'$. 
The following variants of Mader's conjecture have been posed and shown to be true
for $k = 1, 2$ in \cite{H}; more precisely, when $k = 1,2$, 
the minimum degree condition can be weakened to 
$\delta(G) \geq \max\{k+\Delta(T),m-1\}$. 

\begin{conjecture} {\rm (Hasunuma \cite{H})} \label{H1}
For any tree $T$ of order $m$,
every $k$-connected (respectively, $k$-edge-connected) graph $G$ with $\delta(G) \geq k+m -1$ 
contains a subtree $T' \cong T$ 
such that $G-E(T')$ is $k$-connected (respectively, $k$-edge-connected).
\end{conjecture}

Conjecture \ref{H1} generalizes Proposition \ref{tree} and the following two results 
shown in 1969 and 1972 by Halin and Lick, respectively.

\begin{theorem} {\rm (Halin \cite{HR})}
Every $k$-connected graph $G$ with $\delta(G) \geq k+1$ 
contains an edge $uv$ such that $G-uv$ is $k$-connected.
\end{theorem}

\begin{theorem} {\rm (Lick \cite{L})}
Every $k$-edge-connected graph $G$ with $\delta(G) \geq k+1$ 
contains an edge $uv$ such that $G-uv$ is $k$-edge-connected.
\end{theorem}

It has been proved in \cite{H} 
that if two well-known conjectures which strengthen Proposition \ref{tree}, 
namely, the Erd\H{o}s-S\'{o}s conjecture and the Loebl-Koml\'{o}s-S\'{o}s conjecture, are true,
then Conjecture \ref{H1} for $k$-connected graphs is true. 
Conjecture \ref{H1} for $k = 3$ was recently shown to be true in \cite{LLH} and \cite{YT} independently.
For $k \geq 4$, Conjecture \ref{H1} remains open. 

We show that Conjecture \ref{H1} holds for all $k \geq 1$ if we restrict ourselves to cographs.
In particular, we present an improved lower bound on $\delta(G)$ for $k$-edge-connected cographs. 
We also prove the existence of two disjoint connectivity preserving trees when $k = \kappa(G)$.
\begin{itemize}
\item For any trees $T_1$ and $T_2$ of order $m$, every connected cograph $G$
with $\delta(G) \geq \kappa(G)+m-1$ 
contains disjoint subtrees $T'_1 \cong T_1$ and $T'_2 \cong T_2$
such that $\kappa(G-E(T'_1) \cup E(T'_2)) = \kappa(G)$. 
\end{itemize}

We moreover consider the case that $G$ and $G-V(T')$ are $k$-edge-connected
and show the following, where we use Iverson's convention for $[k = 1]$. 

\begin{itemize}
\item For any tree of order $m$, every $k$-edge-connected cograph $G$
with $\delta(G) \geq k+m-[k = 1]$ 
contains a subtree $T' \cong T$ such that $G-V(T')$ is $k$-edge-connected.
\end{itemize}

For $k$-edge-connected cographs where $k \geq 2$, 
this result extends the following result proved in 1986 by Mader. 

\begin{theorem} {\rm (Mader \cite{M0})} \label{M-th}
Every $k$-edge-connected graph $G$ with $\delta(G) \geq k+1$ 
contains a vertex $v$ such that $G-v$ is $k$-edge-connected.
\end{theorem}

Our result and Theorem \ref{M-th} are naturally generalized to the following conjecture.

\begin{conjecture} \label{Con-3}
For any tree $T$ of order $m$,
every $k$-edge-connected graph $G$ with $\delta(G) \geq k+m-[k = 1]$ 
contains a subtree $T' \cong T$ 
such that $G-V(T')$ is $k$-edge-connected.
\end{conjecture}

When $G \cong K_{k+m}$, for any subtree $T'$ of $G$ with $|V(T')| = m$,
$G-V(T') \cong K_k$.
Since $\kappa(K_k) = k-1$, the lower bound of $k+m$ on $\delta(G)$ for $k \geq 2$ 
is tight in general. 
Note that the statement is essentially the same as Mader's conjecture when $k = 1$. 

We furthermore show that every connected cograph is maximally edge-connected
and characterize the super edge-connected cographs.
Based on the characterization, we show the following results on 
a super edge-connectedness keeping tree, where the lower bound of $m+2$ 
on $\delta(G)$ can be shown to be tight.

\begin{itemize}
\item For any tree of order $m$, every super edge-connected cograph $G$
with $\delta(G) \geq m+2$ 
contains a subtree $T' \cong T$ such that $G-V(T')$ is super edge-connected. 
\end{itemize}

This paper is organized as follows.
Section 2 presents structural properties of cographs concerning connectivity,
a characterization of maximally connected cographs
and the proofs of our results on $k$-connected cographs.
A characterization of super edge-connected cographs
and the proofs of our results on $k$-edge-connected cographs
are given in Section 3. 
Section 4 concludes the paper with several remarks.

\section{Mader's Conjecture and Its Variants for $k$-Connected Cographs}

Let $G$ be a cograph.
Since $G$ is $P_4$-free, for any $S \subsetneq V(G)$,
$G-S$ is also $P_4$-free, i.e., $G-S$ is a cograph.
Thus, the class of cographs is closed under vertex-deletion.

\begin{lemma} \label{deletion}
Let $G$ be a cograph.
For any $S \subsetneq V(G)$, $G-S$ is a cograph.
\end{lemma}

By the recursive definition of a cograph $G$,
$G$ can be represented by a binary rooted tree called the {\it parse tree} of $G$ 
in which every leaf is $K_1$
and every non-leaf vertex corresponds to either the union or the join.
Fig. 1 illustrates a cograph $G$ and its parse tree.
Note that a nontrivial cograph $G$ is connected if and only if
the root of the parse tree of $G$ corresponds to the join.
For the parse tree of $G$ and $v \in V(G)$,
by deleting the leaf corresponding to $v$ 
and identifying the parent $u$ of $v$ and another child $w$ of $u$
while keeping the label of $w$,  
the parse tree of $G-v$ is obtained.
From this observation, Lemma \ref{deletion} also follows.

\begin{figure}[h]
\centering
\epsfig{file=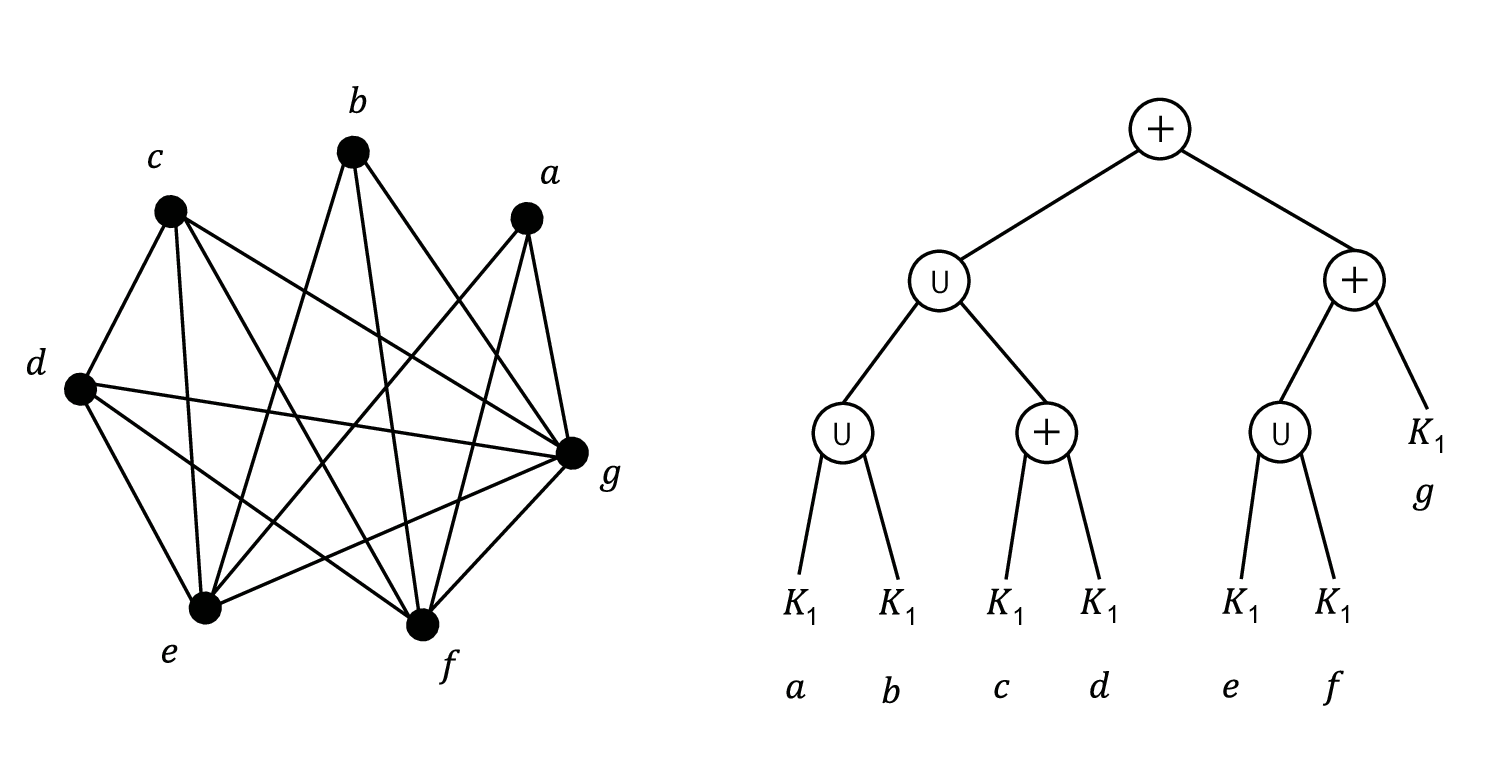,height=55mm}
\caption{A cograph $G$ and its parse tree.}
\end{figure}

Since the union and the join are both associative,
the parse tree of a cograph $G$ can be simplified to the {\it cotree} of $G$ 
by iteratively identifying two vertices $u$ and $v$ if $v$ is a child of $u$
such that both of $u$ and $v$ correspond to the same operation.
That is, the cotree of $G$ is a rooted tree in which every leaf is $K_1$
and every non-leaf vertex is either the union or the join
such that any non-leaf child of the join
(respectively, union) is the union (respectively, join).

Any nontrivial connected (respectively, disconnected) cograph
$G$ can be represented as $G = G_1 + G_2 + \cdots + G_t$
(respectively, $G = G_1 \cup G_2 \cup \cdots \cup G_t$), where $t \geq 2$
and we call each $G_i$ a {\it cocomponent} (respectively, {\it component}) of $G$.
A cocomponent is defined only for a cograph,
while a component is more generally defined as a maximal connected subgraph of a graph. 
Note that a cocomponent is either $K_1$ or a disconnected cograph. 
If a cocomponent is disconnected, then it has at least two components. 
We call a cocomponent with the largest order of a nontrivial connected cograph $G$
a {\it primary cocomponent} of $G$.

Unless stated otherwise, we assume throughout the paper that 
a nontrivial connected cograph $G$ is represented as
$$G = G_1 + G_2 + \cdots + G_{t_G},$$
where 
$$|V(G_1)| \geq |V(G_2)| \geq \cdots \geq |V(G_{t_G})|.$$
Note that we denote by $t_G$ the number of cocomponents of $G$,
i.e., the number of children of the root in the cotree of $G$, 
and $G_1$ is always a primary cocomponent of $G$. 
When $G$ is not a nontrivial connected cograph, 
i.e., $G \cong K_1$ or $G$ is a disconnected cograph,
we consider $G$ itself as a primary cocomponent of $G$ and suppose that $t_G = 1$. 
Let $n'_G$ denote the order of a primary cocomponent of $G$. 
i.e., $n'_G = |V(G_1)|$.
Note that if $n'_G = 1$, then $G$ is a complete graph.
Fig. 2 illustrates the cocomponents $G_1, G_2, G_3$ of the cograph $G$
shown in Fig. 1 and the cotree of $G$, where $n_G = 7$, $n'_G = 4$ and $t_G = 3$.

\begin{figure}[h]
\centering
\epsfig{file=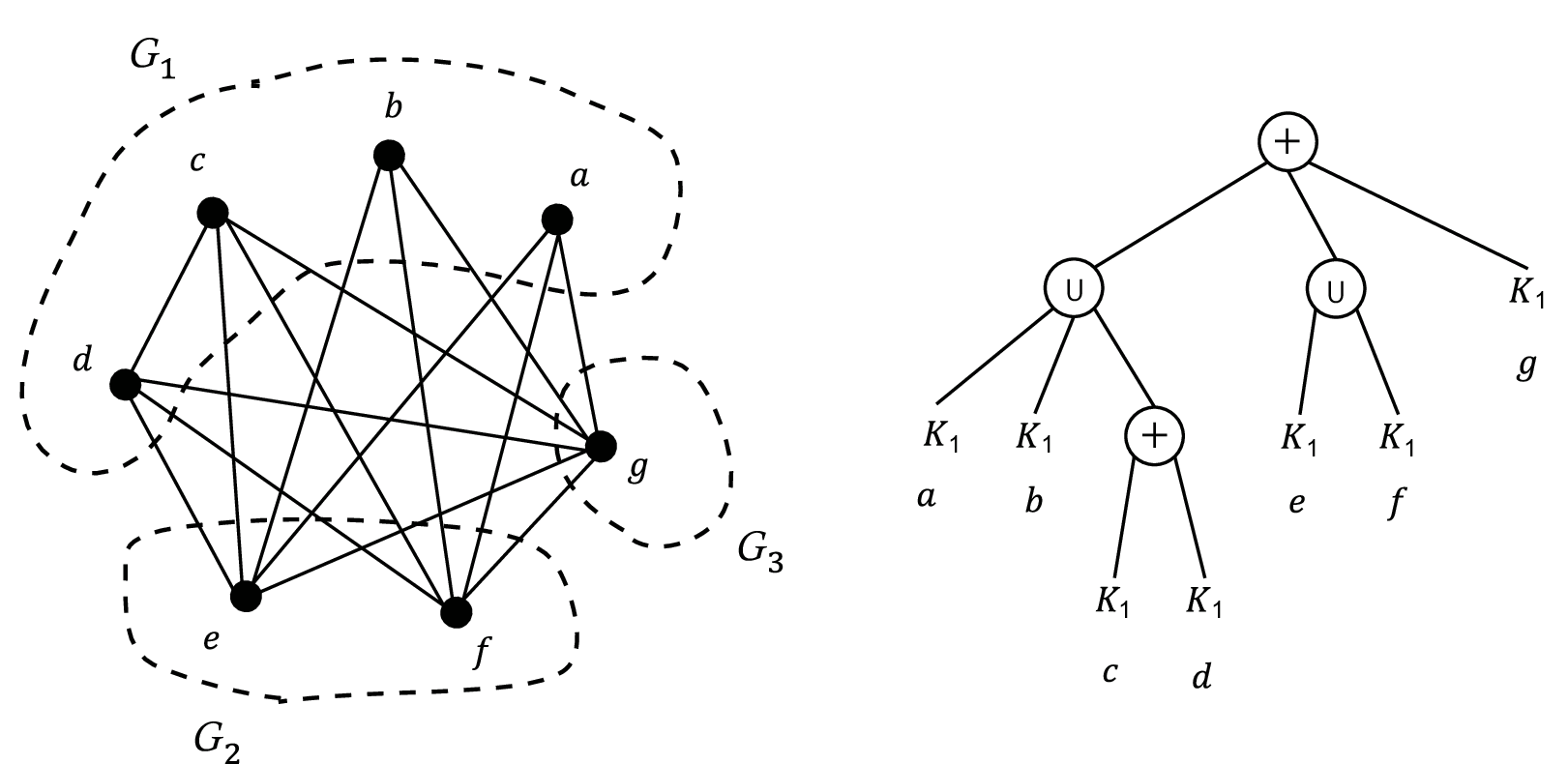,height=55mm}
\caption{The cocomponents $G_1, G_2, G_3$ of the cograph $G$ in Fig. 1 and 
the cotree of $G$, where 
$G_1 = \langle \{a,b,c,d\} \rangle_G$, $G_2 = \langle \{e,f\} \rangle_G$ and $G_3 = \langle \{g\} \rangle_G$. }
\end{figure}

The connectivity of a cograph is determined by its order and
the order of its primary cocomponent. 

\begin{lemma} \label{connectivity}
Let $G$ be a cograph.
Then, $$\kappa(G) = n_G-n'_G.$$
\end{lemma}

\begin{proof}
Let $G$ be a cograph. 
If $G$ is trivial or disconnected, i.e., $t_G = 1$, then $n_G = n'_G$ and the statement holds.
Suppose that $G$ is a nontrivial connected cograph, i.e., $t_G \geq 2$.
For $u \in V(G_i)$ and $v \in V(G_j)$ where $1 \leq i < j \leq t_G$, 
by the definition of the join operation, $uv \in E(G)$ such that 
for any $w \in V(G) \setminus \{u,v\}$, $uw \in E(G)$ or $vw \in E(G)$. 
Thus, for any $S \subsetneq V(G)$ such that 
$V(G_i) \setminus S \neq \emptyset$ and $V(G_j) \setminus S \neq \emptyset$ for $i \neq j$, 
$G-S$ is connected. 
From this observation, for any vertex-cut $X$ of $G$, there exists a cocomponent 
$G_\ell$ such that $V(G) \setminus V(G_\ell) \subseteq X$.
Therefore, for any minimum vertex-cut $Y$ of $G$, there exists a primary cocomponent 
$G_{m}$ such that $Y = V(G) \setminus V(G_{m})$, i.e., $\kappa(G) = n_G-n'_G$. 
\end{proof}

Applying Lemma \ref{connectivity}, we show that 
any graph obtained from a cograph $G$ by deleting 
edges in the cocomponents of $G$ has the same connectivity as $G$.
Note that in contrast to Lemma \ref{deletion}, 
for a cograph $G$ and $F \subset E(G)$, $G-F$ is not always a cograph.

\begin{lemma} \label{edge-deletion}
Let $G$ be a cograph.
For any $F \subseteq \cup_{1 \leq i \leq t_G}E(G_i)$, 
$$\kappa(G-F) = \kappa(G).$$
\end{lemma}

\begin{proof}
Let $G$ be a cograph and $F \subseteq \cup_{1 \leq i \leq t_G}E(G_i)$.
Let $H = G-\cup_{1 \leq i \leq t_G}E(G_i)$.
Since $H \subseteq G-F \subseteq G$ such that $V(H) = V(G-F) = V(G)$, 
it holds that $\kappa(H) \leq \kappa(G-F) \leq \kappa(G)$.
Here, $H$ is the cograph obtained from $G$ by replacing
each cocomponent $G_i$ with $\overline{K_{|V(G_i)|}}$. 
Since $n_H = n_G$ and $n'_H = n'_G$, by Lemma \ref{connectivity}, 
$\kappa(H) = \kappa(G)$.
Therefore, we have $\kappa(G-F) = \kappa(G)$. 
\end{proof}

\begin{corollary} \label{cor-edge-deletion}
Let $G$ be a $k$-connected cograph.
For any $F \subseteq \cup_{1 \leq i \leq t_G}E(G_i)$, 
$G-F$ is $k$-connected.
\end{corollary}

Based on Lemma \ref{connectivity}, we characterize 
the maximally connected cographs as follows. 

\begin{theorem} \label{maximal}
A cograph $G$
is maximally connected
if and only if 
$G$ has a primary cocomponent which has an isolated vertex. 
\end{theorem}

\begin{proof}
Let $G$ be a maximally connected cograph
and $v \in V(G)$ such that ${\rm deg}_G(v)  = \delta(G)$.
Suppose that $v$ is a vertex of a cocomponent $G_\ell$.
Then, it holds that 
$${\rm deg}_G(v) = \sum_{i \neq \ell}|V(G_i)|+{\rm deg}_{G_\ell}(v)
= n_G-|V(G_\ell)|+{\rm deg}_{G_\ell}(v).$$
By Lemma \ref{connectivity}, $\kappa(G) = n_G - n'_G$. 
Since $G$ is maximally connected, $\delta(G) = \kappa(G)$.
Thus, $$|V(G_\ell)| = n'_G + {\rm deg}_{G_\ell}(v).$$
Since $n'_G \geq |V(G_\ell)|$, we have $n'_G = |V(G_\ell)|$ and ${\rm deg}_{G_\ell}(v) = 0$, i.e., 
$G_\ell$ is a primary cocomponent of $G$ and $v$ is an isolated vertex of $G_\ell$.

Conversely, suppose that $G$ has a primary cocomponent which has an isolated vertex $w$. 
For any vertex $v \in V(G_i)$ where $1 \leq i \leq t_G$, 
$${\rm deg}_G(v) = n_G-|V(G_i)|+{\rm deg}_{G_i}(v) \geq n_G - n'_G =  {\rm deg}_G(w).$$
Thus, ${\rm deg}_G(w) = \delta(G)$.
By Lemma \ref{connectivity}, ${\rm deg}_G(w) = \kappa(G)$.
Hence, $G$ is maximally connected.
\end{proof}
\bigskip

From Lemma \ref{connectivity} and Theorem \ref{maximal}, 
the cograph $G$ in Fig. 1 is maximally connected and $\kappa(G) = 3$
(see Fig. 2).  

We next present a sufficient condition for an induced subgraph 
of a $k$-connected cograph to be $k$-connected.

\begin{lemma} \label{k-keeping}
Let $G$ be a $k$-connected cograph. 
Let $S_1 \subseteq V(G_1)$ with $|S_1| \geq \min\{k,|V(G_1)|\}$
and $S_2 \subseteq V(G) \setminus V(G_1)$ with $|S_2| \geq k$.
Then, $\langle S_1 \cup S_2 \rangle_{G}$ is $k$-connected.
\end{lemma}

\begin{proof} 
We consider two possible cases.

\bigskip
\noindent Case 1: $|S_1| \geq k$.

For any $X \subset S_1 \cup S_2$ with $|X| < k$, 
$S_1 \setminus X \neq \emptyset$ and $S_2 \setminus X \neq \emptyset$,
which implies that $\langle S_1 \cup S_2 \rangle_G - X$ is connected.
Therefore, $\langle S_1 \cup S_2 \rangle_G$ is $k$-connected.

\bigskip
\noindent Case 2: $|S_1| = |V(G_1)| < k$.

Let $Y \subset S_1 \cup S_2$ with $|Y| < k$.
Since $S_2 \setminus Y \neq \emptyset$, 
if $S_1 \setminus Y \neq \emptyset$, then $\langle S_1 \cup S_2 \rangle_G - Y$
is connected.
Suppose that $S_1 \setminus Y = \emptyset$, i.e., $S_1 = V(G_1) \subseteq Y$.
Since $$|S_2 \setminus Y| \geq k-(|Y|-|V(G_1)|) = k+|V(G_1)|-|Y| \geq |V(G_1)|+1$$
and $G_1$ is a primary cocomponent,
there is no cocomponent $G_\ell$ with $V(G_\ell) \supseteq S_2 \setminus Y$.
Thus, $S_2 \setminus Y$ has two vertices which are not in the same cocomponent.
This means that 
$\langle S_1 \cup S_2 \rangle_G - Y$ is connected.
Hence, $\langle S_1 \cup S_2 \rangle_G$ is $k$-connected. 
\end{proof}
\bigskip

The following lemma strengthens Proposition \ref{tree}.

\begin{lemma} {\rm (Hasunuma and Ono \cite{HO})} \label{HO-lemma}
Let $T$ be a tree of order $m$ and $T_0$ a subtree of $T$.
If a graph $G$ contains a subtree $T'_0 \cong T_0$ such that ${\rm deg}_G(v) \geq m-1$
for any $v \in V(G) \setminus V(T'_0)$ and any $v \in V(T'_0)$
with ${\rm deg}_{T_0}(\phi^{-1}(v)) < {\rm deg}_T(\phi^{-1}(v))$ where $\phi$ is an isomorphism
from $V(T_0)$ to $V(T'_0)$, then $G$ contains a subtree $T' \cong T$ such that 
$T'_0 \subseteq T'$.
\end{lemma}

From Lemma \ref{HO-lemma}, we have the following weaker but simpler lemma, which will be used to
show our main results.

\begin{lemma} \label{tree2}
Let $G$ be a graph with $\delta(G) \geq m-1$.
Let $T$ be a tree of order $m$ and $T_0$ a subtree of $T$.
If $G$ contains a subtree $T'_0 \cong T_0$, 
then $G$ contains a subtree $T' \cong T$ such that $T'_0 \subseteq T'$.
\end{lemma}

Applying the aforementioned statements, 
we present a minimum degree condition (which contains
all lower bounds on the minimum degree in Conjectures 1,2 and 3) 
for a $k$-connected cograph
to have a connectivity keeping tree.

\begin{theorem} \label{Th1}
For any tree $T$ of order $m$, every $k$-connected cograph $G$ with 
$$\delta(G) \geq 
\left\{ \begin{array}{ll}
k+m-[k = 1] & \mbox{ if } n'_G = 1, \\
k+m-1+\left\lfloor \frac{n'_G}{2} \right\rfloor & \mbox{ if } 2 \leq n'_G \leq k, \\[2mm]
\left\lfloor \frac{3k}{2} \right\rfloor + m-1 & \mbox{ if } k < n'_G < k+m, \\[2mm]
k+m-1 & \mbox{ if } n'_G \geq k+m \\
\end{array} \right.$$
contains a subtree $T' \cong T$ such that $G-V(T')$ is $k$-connected.
\end{theorem}

\begin{proof}
Let $T$ be any tree of order $m$ 
and $G$ a $k$-connected cograph satisfying the minimum degree condition
in the statement. 
Let $S = V(G) \setminus V(G_1)$.
By Lemma \ref{connectivity}, we have $\kappa(G) = |S| \geq k$.

We consider four cases depending on $n'_G = |V(G_1)|$. 

\bigskip
\noindent Case 1: $|V(G_1)| = 1$.

In this case, $G = \underbrace{K_1 + K_1 + \cdots + K_1}_{n_G} \cong K_{n_G}$.
Since $\delta(G) \geq k+m-[k = 1]$,
$G$ contains a subtree $T' \cong T$ such that $G-V(T') \cong K_{n_G-m}$,
where $n_G-m \geq k+1-[k = 1]$.
If $k \geq 2$, then $\kappa(G-V(T')) \geq k$.
When $k = 1$, $G-V(T')$ is 1-connected even if $n_G-m = 1$.
Note that $K_1$ is assumed to be 1-connected.
Thus, $G-V(T')$ is $k$-connected for any $k \geq 1$.

\bigskip
\noindent Case 2: $2 \leq |V(G_1)| \leq k$.

Since $G_1$ is disconnected,
$G_1$ has at least two components.
Let $G'_1$ be a component of $G_1$ with the least order.
Then, $|V(G'_1)| \leq \left\lfloor \frac{|V(G_1)|}{2} \right\rfloor$ and 
$\Delta(G'_1) \leq \left\lfloor \frac{|V(G_1)|}{2} \right\rfloor -1$.
Since $\delta(G) \geq k+m-1+ \left\lfloor \frac{|V(G_1)|}{2} \right\rfloor$,
any vertex of $G'_1$ is adjacent to at least $k+m$ vertices in $S$, which means that
$|S| \geq k+m$.
Since $|V(G_1)| \leq k$, 
$\delta(\langle S \rangle_G) \geq m-1+ \left\lfloor \frac{|V(G_1)|}{2} \right\rfloor$. 
Thus, 
$\langle S \rangle_G$ contains a subtree $T' \cong T$.
Here, it holds that $|S \setminus V(T')| \geq k$ and 
$G-V(T') = \langle V(G_1) \cup (S \setminus V(T')) \rangle_G$.
Therefore, it follows from Lemma \ref{k-keeping} that $G-V(T')$ is $k$-connected.

\bigskip
\noindent Case 3: $k < |V(G_1)| < k+m$.

Let $|V(G_1)|= k+p$ where $1 \leq p \leq m-1$ 
and $G'_1$ a component of $G_1$ with the least order.
Then, $|V(G'_1)| \leq \left\lfloor \frac{k+p}{2} \right\rfloor$.
Since $\delta(G) \geq \left\lfloor \frac{3k}{2} \right\rfloor +m-1$, 
$$|S| \geq \left(\left\lfloor \frac{3k}{2} \right\rfloor +m-1\right) - 
\left(\left\lfloor \frac{k+p}{2} \right\rfloor-1\right) 
\geq \frac{3k}{2}-\frac{1}{2}+m-\frac{k+p}{2} \geq k+m-p.$$
Let $T_0$ be a subtree of $T$ with $|V(T_0)| = m-p$. 
Since $\delta(\langle S \rangle_G) \geq \delta(G)-|V(G_1)| \geq \lfloor \frac{k}{2} \rfloor + m-1-p$, 
$\langle S \rangle_G$ contains a subtree $T'_0 \cong T_0$.
Let $S' \subset S \setminus V(T'_0)$ with $|S'| = k$.
Then, $\delta(G-S') \geq \lfloor \frac{k}{2} \rfloor +m-1$.
Thus, by Lemma \ref{tree2}, 
we can extend $T'_0$ to $T' \cong T$ in $G-S'$.
Since $|V(G_1) \cap V(T')| \leq p$, we have $|V(G_1) \setminus V(T')| \geq k$.
Moreover, $|S \setminus V(T')| \geq |S'| \geq k$.
Hence, from Lemma \ref{k-keeping}, $G-V(T')$ is $k$-connected.

\bigskip
\noindent Case 4: $|V(G_1)| \geq k+m$.

Let $S' \subseteq S$ with $|S'| = k$.
Since $\delta(G) \geq k+m-1$, $\delta(G-S') \geq m-1$. 
Thus, $G-S'$ contains a subtree $T' \cong T$.
Since $|V(G_1) \setminus V(T')| \geq k$ and $|S \setminus V(T')| \geq k$, 
by Lemma \ref{k-keeping}, $G-V(T')$ is $k$-connected.
\end{proof}
\bigskip

For any $k \geq 1$, the maximum of the four lower bounds on $\delta(G)$ in Theorem \ref{Th1} 
is $\left\lfloor \frac{3k}{2} \right\rfloor +m-1$. 
Thus, Mader's conjecture holds for $k$-connected cographs.

\begin{corollary} \label{Co1}
For any tree $T$ of order $m$, every $k$-connected cograph $G$
with $\delta(G) \geq \left\lfloor \frac{3k}{2} \right\rfloor + m-1$
contains a subtree $T' \cong T$ such that $G-V(T')$ is $k$-connected.
\end{corollary}

The condition $n'_G > \delta(G)$ is equivalent to $n_G > \kappa(G)+\delta(G)$
since $\kappa(G) = n_G-n'_G$.
Thus, from the fourth condition on $\delta(G)$ in Theorem \ref{Th1}, 
we have the following corollary.

\begin{corollary} \label{Co2}
For any tree $T$ of order $m$, every $k$-connected cograph $G$
with $n_G > \kappa(G)+\delta(G)$ and $\delta(G) \geq k+m-1$
contains a subtree $T' \cong T$ such that $G-V(T')$ is $k$-connected. 
\end{corollary}

Any graph with $\delta(G) \geq \frac{n_G}{2}$ where $n_G \geq 3$ is called a Dirac graph
named after Diarc's famous sufficient condition for a graph to be Hamiltonian \cite{D}. 
Concerning a connectivity keeping tree, it has been shown in \cite{H2} that 
every $k$-connected Dirac graph of large order has a connectivity preserving Hamiltonian
cycle.
In contrast to this result, the following result is obtained from Corollary \ref{Co2}
by applying the fact that $\kappa(G) \leq \delta(G)$.

\begin{corollary} \label{CoD}
For any tree $T$ of order $m$, every $k$-connected non-Dirac cograph $G$
with $\delta(G) \geq k+m-1$ 
contains a subtree $T' \cong T$ such that $G-V(T')$ is $k$-connected. 
\end{corollary}

For general $k$-connected graphs, the lower bound of $\left\lfloor \frac{3k}{2} \right\rfloor +m-1$ 
on $\delta(G)$
in Mader's conjecture is certainly best possible;
however, the $k$-connected graphs shown in \cite{M1} for the tightness are not cographs. 
We then discuss tightness of the lower bounds on $\delta(G)$
in Theorem \ref{Th1}.

\begin{proposition}
In Theorem \ref{Th1}, 
\begin{enumerate}
\item the lower bound of $k+m-[k = 1]$ on $\delta(G)$ for $n'_G = 1$ is tight, 
\item the lower bound of $\left\lfloor \frac{3k}{2} \right\rfloor+m-1$ on $\delta(G)$ for 
$n'_G = k$ is at most $\lfloor \frac{m}{2} \rfloor$ greater than the optimum, 
\item the lower bound of $\left\lfloor \frac{3k}{2} \right\rfloor +m-1$ on $\delta(G)$ for
$k < n'_G < k+m$ is at most $\frac{m-1}{4}$ greater than the optimum, 
\item the lower bound of $k+m-1$ on $\delta(G)$ for $n'_G \geq k+m$ is tight. 
\end{enumerate}
\end{proposition}

\begin{proof}
Theorem \ref{Th1} is a statement for any $m \geq 1$ and any $k \geq 1$.
Thus, if there exists a $k$-connected cograph $H$ for specific values $k$ and $m$
such that for any subtree $T' \subset H$ with $|V(T')| = m$, $H-V(T')$ is not $k$-connected
(or $H-V(T')$ is not well-defined), 
then the optimum lower bound on the minimum degree is at least $\delta(H)+1$.

\bigskip
\noindent Case 1: $n'_G = 1$.

Let $H = K_{k+m-[k = 1]}$.
Then, $n'_H = 1$ and $\kappa(H) = \delta(H) = k+m-1-[k = 1]$.
For any subtree $T' \subset H$ with $|V(T')| = m$, 
if $k \geq 2$ (respectively, $k = 1$), then $H-V(T') \cong K_{k}$ 
is not $k$-connected (respectively, $H-V(T')$ is not well-defined).
Therefore, the lower bound of $k+m-[k = 1]$ is tight. 

\bigskip
\noindent Case 2: $n'_G = k$. 

Suppose that $m < \lfloor \frac{k}{2} \rfloor$ and $k \geq 4$. 
Let $H = H_1 + H_2 + \underbrace{K_1+K_1+\cdots+K_1}_{\lceil \frac{m}{2} \rceil -1}$, where 
$$H_1 \cong H_2 \cong K_{\lfloor \frac{k}{2} \rfloor} \cup K_{\lceil \frac{k}{2} \rceil}.$$
Then, $n'_H = k$, $\kappa(H) = k+\lceil \frac{m}{2} \rceil -1$ and 
$\delta(H) = \left\lfloor \frac{3k}{2} \right\rfloor + \lceil \frac{m}{2} \rceil -2$.
Let $T' \subset H$ be a subtree of order $m$ such that 
$|V(T') \cap V(H_1)| = m_1$ and $|V(T') \cap V(H_2)| = m_2$. 
Since $m < \lfloor \frac{k}{2} \rfloor$, both of $H_1 - V(T') \cap V(H_1)$ 
and $H_2 - V(T') \cap V(H_2)$ are disconnected, i.e.,
they are cocomponents of the cograph $H-V(T')$.
Thus, the order of a primary cocomponent of $H-V(T')$ is $\max\{k-m_1, k-m_2\}$. 
Therefore, 
$$\begin{array}{ll}
\kappa(H-V(T')) & = n_{H-V(T')}-n'_{H-V(T')} \\
& = (2k+\lceil \frac{m}{2} \rceil-1-m)-\max\{k-m_1, k-m_2\} \\
& =k-1-\lfloor \frac{m}{2} \rfloor +\min\{m_1,m_2\} \\
& < k. \\
\end{array}$$
Thus, the optimum lower bound is at least 
$\left\lfloor \frac{3k}{2} \right\rfloor + \lceil \frac{m}{2} \rceil -1$.

\bigskip
\noindent Case 3: $k < n'_G < k+m$. 

Suppose that $5 \leq m \leq 2k$, $m \bmod{4} = 1$ and $k$ is even.
Let $H = H_1+H_2$, where 
$$H_1 \cong H_2 \cong 
K_{\frac{k}{2} + \frac{m-1}{4}} 
\cup K_{\frac{k}{2} + \frac{m-1}{4}}.$$
Then, $n'_H = \kappa(H) = k+\frac{m-1}{2}$ and 
$\delta(H) = \frac{3k}{2} +\frac{3(m-1)}{4}-1$.
Let $T' \subset H$ be a subtree of order $m$.
We may assume, without loss of generality, that
$|V(H_1) \cap V(T')| \geq \lceil \frac{m}{2} \rceil$
and $|V(H_2) \cap V(T')| \leq \lfloor \frac{m}{2} \rfloor$. 
Since $|V(H_1)| = k+\frac{m-1}{2} < k+ \lceil \frac{m}{2} \rceil$, 
$|V(H_1) \setminus V(T')|  < k$.
Thus, if $\langle V(H_2) \setminus V(T') \rangle_{H}$ is disconnected,
then $\kappa(H-V(T')) < k$.
Assume that $\langle V(H_2) \setminus V(T') \rangle_{H}$ is connected.
Then, $|V(H_2) \setminus V(T')| \leq \frac{k}{2} + \frac{m-1}{4}$,
i.e., $|V(H_2) \cap V(T')| \geq \frac{k}{2} + \frac{m-1}{4}$.
Since $|V(H_2) \cap V(T')| \leq \lfloor \frac{m}{2} \rfloor$, we have
$\frac{k}{2} + \frac{m-1}{4} \leq \lfloor \frac{m}{2} \rfloor$, i.e., $2k \leq m-1$, 
which contradicts the assumption that $m \leq 2k$.
Therefore, we have $\kappa(H-V(T')) < k$.
Thus, the optimum lower bound is at least 
$\frac{3k}{2} +\frac{3(m-1)}{4}$.

\bigskip
\noindent Case 4: $n'_G \geq k+m$. 

Suppose that $k \leq 2m-3$. 
Let $H = H_1+\underbrace{K_1+K_1+\cdots+K_1}_{k}$, where  
$$H_1 = K_{m-1} \cup K_{m-1} \cup K_{m-1}.$$
Then, $n'_H = 3(m-1) \geq k+m$, $\kappa(H) = k$ and 
$\delta(H) = k+m-2$.
For any subtree $T' \subset H$ with $|V(T')| = m$, 
$(V(H) \setminus V(H_1)) \cap V(T') \neq \emptyset$
and $H_1 - V(T') \cap V(H_1)$ is disconnected. 
Thus, $H_1 - V(T') \cap V(H_1)$ is a primary cocomponent of
$H-V(T')$ and $\kappa(H-V(T')) < k$.
Hence, the lower bound of $k+m-1$ is tight.
\end{proof}
\bigskip

We next show that under a slightly stronger condition than that of Mader's conjecture, 
there are two disjoint connectivity keeping trees when $k = \kappa(G)$. 

\begin{theorem} \label{Th2}
For any trees $T_1$ and $T_2$ of order $m$,
every connected cograph $G$ with 
$\delta(G) \geq \left\lceil \frac{3\kappa(G)}{2} \right\rceil +m-1$ 
contains $T'_1 \cong T_1$ and $T'_2 \cong T_2$ with 
$V(T'_1) \cap V(T'_2) = \emptyset$ such that 
$\kappa(G-V(T'_1) \cup V(T'_2)) = \kappa(G)$. 
\end{theorem}

\begin{proof}
Let $G$ be a connected cograph $G$ with 
$\delta(G) \geq \left\lceil \frac{3\kappa(G)}{2} \right\rceil +m-1$. 
Let $T_1$ and $T_2$ be any trees of order $m$. 
If $n'_G = 1$, then $G \cong K_{n_G}$ and 
$\kappa(G) = \delta(G) \geq \left\lceil \frac{3\kappa(G)}{2} \right\rceil +m-1$,
which is a contradiction.
Thus, $n'_G \geq 2$, i.e., $G_1$ is disconnected and has at least two components.
Now let $G'_1$ and $G''_1$ be components of $G_1$. 
Since both $G'_1$ and $G''_1$ have minimum degree at least 
$\delta(G)-\kappa(G) \geq \left\lceil \frac{\kappa(G)}{2} \right\rceil +m-1$, 
$G'_1$ and $G''_1$ contain subtrees $T'_1 \cong T_1$ and $T'_2 \cong T_2$, respectively,
such that 
$|V(G'_1) \setminus V(T'_1)| \geq \left\lceil \frac{\kappa(G)}{2} \right\rceil$
and $|V(G''_1) \setminus V(T'_2)| \geq \left\lceil \frac{\kappa(G)}{2} \right\rceil$. 
Thus, 
$G_1 - V(T'_1) \cup V(T'_2)$ is a disconnected cograph with at least $\kappa(G)$ vertices.
Therefore, 
$G_1 - V(T'_1) \cup V(T'_2)$ is a primary cocomponent of 
$G-V(T'_1) \cup V(T'_2)$.
Hence, $\kappa(G-V(T'_1) \cup V(T'_2)) = \kappa(G)$.
\end{proof}
\bigskip

Suppose that $k \geq 3$, $k$ is odd and $m < \lfloor \frac{k}{2} \rfloor$.
Let $H = H_1 + H_2$, where 
$$H_1 \cong K_{\lfloor \frac{k}{2} \rfloor +m} \cup K_{\lfloor \frac{k}{2} \rfloor +m}, \ \ \ 
H_2 \cong K_{\lfloor \frac{k}{2} \rfloor} \cup K_{\lceil \frac{k}{2} \rceil}.$$
Then, $\kappa(H) = k$ and $\delta(H) = \left\lfloor \frac{3k}{2} \right\rfloor +m-1$. 
Let $T'_1$ and $T'_2$ be disjoint trees of order $m$ in $H$.
Since $2m < \lfloor \frac{k}{2} \rfloor + m$, 
$H_1 - (V(T'_1) \cup V(T'_2)) \cap V(H_1)$ is disconnected. 
If $V(H_2) \cap (V(T'_1) \cup V(T'_2)) \neq \emptyset$,
then $H_1-(V(T'_1) \cup V(T'_2)) \cap V(H_1)$ has at least $2 \lfloor \frac{k}{2} \rfloor+1 = k$ 
vertices and it is a primary cocomponent of $H-V(T'_1) \cup V(T'_2)$. 
If $V(H_2) \cap (V(T'_1) \cup V(T'_2)) = \emptyset$,
then $|V(H_1) \setminus (V(T'_1) \cup V(T'_2))| = 2 \lfloor \frac{k}{2} \rfloor = k-1$ and 
$H_2$ is a primary cocomponent of $H-V(T'_1) \cup V(T'_2)$.
Thus, in either case, we have $\kappa(H-V(T'_1) \cup V(T'_2)) < k$.
Therefore, the lower bound on $\delta(G)$ in Theorem \ref{Th2} is tight.

Any nontrivial tree $T$ is a bipartite graph, i.e.,  
$V(T)$ can be divided into $V_1$ and $V_2$ so that 
any edge of $T$ joins a vertex in $V_1$ and a vertex in $V_2$.
This partition $(V_1,V_2)$ is called the {\it bipartition} of $T$.
Let $K_{n_1,n_2,\ldots,n_r} = \overline{K_{n_1}} + \overline{K_{n_2}} + \cdots + \overline{K_{n_r}}$. 
That is, $K_{n_1,n_2, \ldots, n_r}$ is the complete $r$-partite graph with 
partite sets of cardinalities $n_1,n_2,\ldots,n_r$.  
For any tree $T$ with bipartition $(V_1,V_2)$, 
$T \cong T' \subseteq K_{|V_1|,|V_2|}$.

In what follows, we show that under the condition that $\delta(G) \geq m-1$,
every maximally connected cograph $G$ has a maximal connectedness keeping tree, 
except for the trivially inappropriate case that $G \cong K_{m}$. 
The lower bound of $m-1$ on $\delta(G)$ is clearly tight. 

\begin{theorem} \label{max-con-tree}
For any tree $T$ of order $m$, 
every maximally connected cograph $G \not\cong K_m$ with 
$\delta(G) \geq m-1$
contains a subtree $T' \cong T$ such that $G-V(T')$ is maximally connected.
\end{theorem}

\begin{proof}
Let $T$ be a tree of order $m$. 
Let $G \not\cong K_m$ be a maximally connected cograph with $\delta(G) = k+m-1$, 
where $k$ is a nonnegative integer.
By Theorem \ref{maximal}, a primary cocomponent of $G$ has an isolated vertex.
Without loss of generality, we may assume that $G_1$ has an isolated vertex $x$.
Let $|V(G_1)| = k+p$, where $p$ may be a negative integer.
Let $S = V(G) \setminus V(G_1)$.
Since $G$ is maximally connected, $|S| = \kappa(G) = \delta(G) = k+m-1$.

Suppose that $m = 1$, i.e., $T \cong K_1$. 
Since $G \not\cong K_1$, $n_G \geq 2$.
If $\kappa(G) = 0$, i.e., $S = \emptyset$, then for any $y \in V(G) \setminus \{x\}$,
$\kappa(G-y) = \delta(G-y) = 0$.
If $\kappa(G) \geq 1$, then 
for any $z \in S$, 
$\kappa(G-z) = \delta(G-z) = \kappa(G)-1$.
Hence, the statement holds when $m = 1$.

Suppose that $m \geq 2$.
Let $(V_1,V_2)$ be the bipartition of $T$, where $|V_1| \geq |V_2|$.
We consider three cases depending on $p$.

\bigskip
\noindent Case 1: $p \geq m-1$.

Let $S_1 \subseteq S$ with $|S_1| = |V_1|$ and 
$S_2 \subseteq V(G_1) \setminus \{x\}$ with $|S_2| = |V_2|$.
Note that $|V(G_1)| > |V_2|$ and $S_2$ is well-defined; 
for otherwise $k+m-1 \leq k+p = |V(G_1)| \leq |V_2| \leq \frac{m}{2}$, i.e., $m = 2$,
$k = 0$ and $p = 1$, which contradicts the assumption that $G \not\cong K_m$.
Since $K_{|V_1|,|V_2|} \subseteq \langle S_1 \cup S_2 \rangle_G$, 
$T \cong T' \subseteq \langle S_1 \cup S_2 \rangle_G$, where
$V(T') = S_1 \cup S_2$.
Note that $G_1-S_2$ is a cocomponent of $G-V(T')$. 
Since $|V(G_1)| \geq |S|$ and $|S_1| \geq |S_2|$, $|V(G_1) \setminus S_2| \geq |S \setminus S_1|$.
Therefore, $G_1-S_2$ is a primary cocomponent of $G-V(T')$ such that $x \in V(G_1-S_2)$.
Hence, $G-V(T')$ is maximally connected.

\bigskip
\noindent Case 2: $p \leq 0$.

In this case, it holds that $\delta(\langle S \rangle_G) \geq k+m-1-|V(G_1)| \geq m-1$.
Thus, $T \cong T' \subseteq \langle S \rangle_G$ and 
$G_1$ is still a primary cocomponent of $G-V(T')$.
Therefore, $G-V(T')$ is maximally connected.

\bigskip
\noindent Case 3: $1 \leq p \leq m-2$.

Let $q = |V_1|-|V_2| \geq 0$.
Since $|V_1| + |V_2| = m$, $|V_1| = \frac{m+q}{2}$ and $|V_2| = \frac{m-q}{2}$.
Note that $m$ and $q$ have the same parity.
Two cases are now considered, depending on $q$. 

\bigskip
\noindent Case 3.1: $q \geq m-p-1$.

In this case, $|V_1|-|V_2| \geq |S| - |V(G_1)|$, i.e., $|V(G_1)|-|V_2| \geq |S|-|V_1| \geq k \geq 0$.
Note that $|V(G_1)| > |V_2|$, since otherwise
$|V(G_1)| = |V_2|$ and $|S| = |V_1|$, i.e, $|V(G)| = m$ and $G \cong K_m$, which is a contradiction.
Similarly to Case 1, by letting $S_1 \subseteq S$ with $|S_1| = |V_1|$ and 
$S_2 \subseteq V(G_1) \setminus \{x\}$ with $|S_2| = |V_2|$,
we have $T \cong T' \subseteq \langle S_1 \cup S_2 \rangle_G$ and 
$|V(G_1) \setminus S_2| \geq |S \setminus S_1|$.
Thus, $T'$ is a desired tree.

\bigskip
\noindent Case 3.2: $q < m-p-1$.

Let $r = \left\lceil \frac{m-p-q-1}{2} \right\rceil > 0$.
Note that 
\begin{equation} \label{m-1}
m-1 = p+q+2r -[p \mbox{ even}],
\end{equation}
where we use Iverson's convention for $[p \mbox{ even}]$, i.e., 
$[p \mbox{ even}] = 1$ (respectively, 0) if $p$ is even (respectively, odd) and
\begin{equation} \label{V1V2}
|V_1| = q+r +\left\lfloor \frac{p+1}{2} \right\rfloor, \ \ \ 
|V_2| = r + \left\lfloor \frac{p+1}{2} \right\rfloor.
\end{equation}
Since $|S| = k+m-1$ and $|V(G_2)| \leq |V(G_1)| = k+p$,
we have 
$$|S \setminus V(G_2)| \geq (k+m-1)-(k+p) = m-1-p \geq q+2r-1 > 0.$$
Thus, $t_G \geq 3$, i.e., the cocomponent $G_3$ exists in $G$.
Fig. 3 illustrates the cardinalities of $V_1, V_2, V(G_1)$ and $S$ 
with the values $k,m,p,q$ and $r$. 
We distinguish two cases.

\begin{figure}[t]
\centering
\epsfig{file=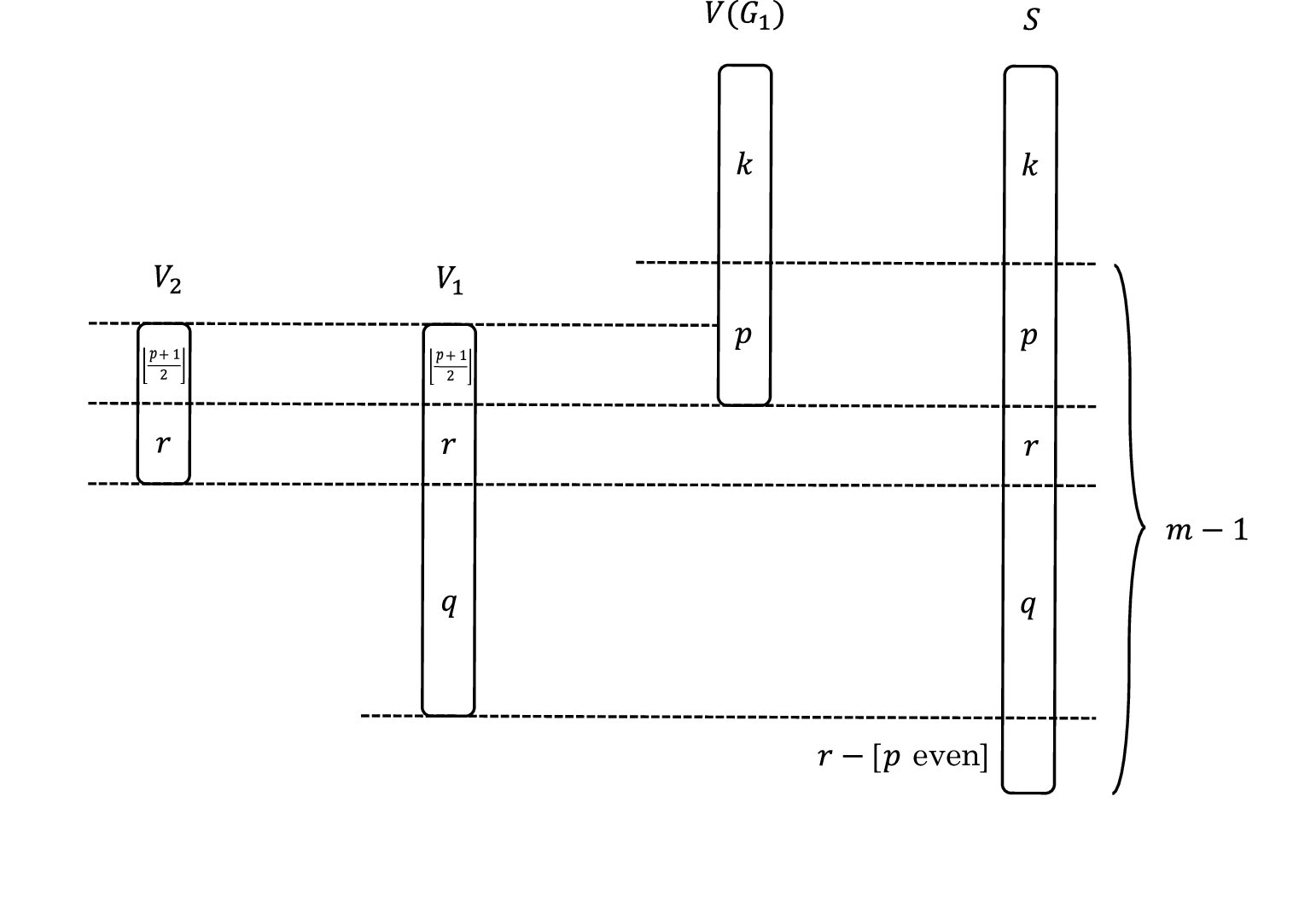,height=100mm}
\caption{The cardinalities of $V_1,V_2, V(G_1)$ and $S$ with the values $k,m,p,q$ and $r$ in 
Case 3.2 in the proof of Theorem \ref{max-con-tree}.}
\end{figure}

\bigskip
\noindent Case 3.2.1: $|V(G_2)| \geq |V_2|$.

Let $S_2 \subseteq V(G_2)$ with $|S_2| = |V_2|$,
$S_{1,1} \subseteq S \setminus V(G_2)$ with $|S_{1,1}| = q+r$
and $S_{1,2} \subseteq V(G_1) \setminus \{x\}$ with $|S_{1,2}| = \lfloor \frac{p+1}{2} \rfloor$.
Note that $|V(G_1)| > \frac{p+1}{2}$; for otherwise 
$k+p \leq \frac{p+1}{2}$, i.e., $p = 1$ and $k = 0$, which is contradictory.
Since $G = G_1 + G_2 + \langle S \setminus V(G_2) \rangle_G$, 
$$\langle S_{1,2} \cup S_2 \cup S_{1,1} \rangle_G 
\supseteq K_{\lfloor \frac{p+1}{2} \rfloor, |V_2|, q+r}
\supseteq K_{|V_1|,|V_2|}.$$
Thus, $\langle S_{1,2} \cup S_2 \cup S_{1,1} \rangle_G$ contains a subtree $T' \cong T$,
where $V(T') = S_{1,2} \cup S_2 \cup S_{1,1}$. 
Let $G'_1 = G_1 - S_{1,2}$ and $S' = S \setminus (S_2 \cup S_{1,1})$. 
Then, from (\ref{m-1}), it follows that
$$|S'| = k+m-1-\left(r+\left\lfloor \frac{p+1}{2} \right\rfloor + q+r\right) 
\leq k+p-\left\lfloor \frac{p+1}{2} \right\rfloor = |V(G'_1)|.$$
Thus, $G'_1$ is a primary cocomponent of $G-V(T')$ with the isolated vertex $x$.
Therefore, we have the desired result.

\bigskip
\noindent Case 3.2.2: $|V(G_2)| < |V_2|$.

Since $|V(G_2)| \geq |V(G_3)|$, from (\ref{V1V2}), we have 
$$|S|-(|V(G_2)|+|V(G_3)|) \geq k+m-1-2(|V_2|-1) = k+q+1 > 0.$$
Thus, we have $t_G \geq 4$.
If $|V(G_2)|+|V(G_4)| < |V_2|$, then we similarly have $t_G \geq 5$ and 
$|V(G_3)|+|V(G_5)| \leq |V(G_2)|+|V(G_4)| < |V_2|$, 
which also implies that $t_G \geq 6$.
From these observations, we may suppose that
there exists an integer $2 \leq s \leq \frac{t_G}{2}$ such that 
$$\sum_{1 \leq i <  s}|V(G_{2i})| < |V_2|,\ \ \ \sum_{1 \leq i \leq  s}|V(G_{2i})| \geq |V_2|.$$
Now let $S_{2} = \cup_{1 \leq i \leq s}V(G_{2i})$  
and $S_{1} = S \setminus S_2$.
If $|S_1| \geq |V_1|$, then  
$$\langle S \rangle_G \supset 
K_{|V(G_2)|,|V(G_{3})|,\ldots,|V(G_{t_G})|} 
\supset K_{|S_1|,|S_2|} \supset K_{|V_1|,|V_2|} \supseteq T' \cong T,$$
and $G_1$ is still a primary cocomponent of $G-V(T')$; hence $G-V(T')$ is maximally connected.

Suppose that $|S_1| < |V_1|$.
Let $\ell_1 = |V_1|-|S_1| > 0$ and $\ell_2 = |V_2|-\sum_{1 \leq i < s}|V(G_{2i})| > 0$.
Note that $\ell_1+\ell_2 = |V(G_{2s})|-k+1$ and $\max\{\ell_1,\ell_2\} \leq |V(G_{2s})|$.
We finally consider two cases depending on $\ell_1$ and $\ell_2$.

\bigskip
\noindent Case 3.2.2.1: $\ell_2 \geq \ell_1$.

Let $L_2 \subseteq V(G_{2s})$ with $|L_2| = \ell_2$ and $S'_2 = \cup_{1 \leq i < s}V(G_{2i}) \cup L_2$.
Also, let $U_1 \subseteq V(G_1) \setminus \{x\}$ with 
$|U_1| = \ell_1$.
Note that 
$$|V(G_1)|-|U_1| \geq |V(G_{2s})|-\ell_1 = k+\ell_2-1.$$
If $k > 0$ or $\ell_2 \geq 2$, then $U_1$ is well-defined.
Suppose that $k = 0$ and $\ell_2 = 1$.
Then, $|U_1| = \ell_1 = 1$.
If $|V(G_1)| = 1$, then $n_G = m$ which is a contradiction.
Thus, $|V(G_1)| \geq 2$ and $U_1$ is well-defined.  
Since $|S'_2| = |V_2|$ and $|S_1 \cup U_1| = |V_1|$, 
we have 
$$\langle U_1 \cup S_1 \cup S'_2 \rangle_G \supseteq K_{|U_1|,|S_1|,|S'_2|} \supseteq K_{|V_1|,|V_2|}
\supseteq T' \cong T,$$
where $V(T') = U_1 \cup S_1 \cup S'_2$.
Since 
$$|V(G_1) \setminus U_1| \geq  |V(G_{2s})|-\ell_1 
\geq  |V(G_{2s})|-\ell_2 = |S \setminus  (S_1 \cup S'_2)|,$$ 
$G_1 - U_1$ is a primary cocomponent of $G-V(T')$.
Therefore, $T'$ is a desired tree.

\bigskip
\noindent Case 3.2.2.2: $\ell_2 < \ell_1$.

Let $L_1 \subseteq V(G_{2s})$ with $|L_1| = \ell_1$ and 
$S'_1 = S_1 \cup L_1$.
Also, let $S'_2 = \cup_{1 \leq i < s}V(G_{2i})$
and $U_2 \subseteq V(G_1) \setminus \{x\}$ with $|U_2| = \ell_2$.
Similarly to Case 3.2.2.1, 
we have 
$\langle S'_1 \cup U_2 \cup S'_2 \rangle_G \supset K_{|V_1|,|V_2|} \supseteq T' \cong T$
where $V(T') = S'_1 \cup U_2 \cup S'_2$
such that $|V(G_1) \setminus U_2| > |S \setminus  (S'_1 \cup S'_2)|$.
Thus, $G_1 - U_2$ is a primary cocomponent of $G-V(T')$.
Hence, we have the desired result.
\end{proof}
\bigskip

By Whitney's famous theorem which characterizes $k$-connected graphs \cite{W}, 
a maximally connected graph can be characterized
as a graph $G$ such that 
for any two vertices in $G$, there are $\delta(G)$
internally disjoint paths between them.
Concerning the maximally connected graphs,
a much stronger notion is known.
If for any two vertices $u, v$ in $G$, there are $\min\{ {\rm deg}_G(u), {\rm deg}_G(v)\}$
internally disjoint paths between $u$ and $v$, then
$G$ is called {\it ideally connected}.
Thus, the class of ideally connected graphs is properly included in the class of
maximally connected graphs. 

Jordaan recently characterized the ideally connected cographs as follows,
where $2K_2 = K_2 \cup K_2$.

\begin{theorem} {\rm (Jordaan \cite{J})} \label{ideal-tree}
A cograph $G$ is ideally connected if and only if $G$ is $2K_2$-free.
\end{theorem}

From Theorem \ref{ideal-tree}, 
the following result on an ideal connectedness keeping tree
for cographs is immediately obtained.

\begin{corollary}
For any tree $T$ of order $m$, 
every ideally connected cograph $G \not\cong K_m$ with 
$\delta(G) \geq m-1$
contains a subtree $T' \cong T$ such that $G-V(T')$ is ideally connected.
\end{corollary}

Note that any subtree of an ideally connected cograph is actually 
an ideal connectedness keeping tree.

At the last of this section,
we consider Conjecture \ref{H1} for $k$-connected cographs. 
We first show the following.
The tightness of the lower bound of $k+m-1$ on $\delta(G)$ 
follows from the fact that
$\delta(G-V(T')) = k-1$ for any subtree $T' \cong T$
when $G \cong K_{k+m-1}$ and $T$ is the star of order $m$, i.e., $\Delta(T) = m-1$. 

\begin{theorem} \label{Th3}
For any tree $T$ of order $m$, every $k$-connected cograph $G$
with $\delta(G) \geq k+m-1$ 
contains a subtree $T' \cong T$ such that $G-E(T')$ is $k$-connected. 
\end{theorem}

\begin{proof}
Let $G$ be a $k$-connected cograph with $\delta(G) \geq k+m-1$.
Let $S = V(G) \setminus V(G_1)$ and $S' \subseteq S$ with $|S'| = k$.
Since $\delta(G-S') \geq m-1$, $G-S'$ contains a subtree $T' \cong T$.
By Corollary \ref{cor-edge-deletion}, $G-E(G_1) \cap E(T')$ is $k$-connected.
Since $|S \setminus V(T')| \geq |S'| \geq k$, from Lemma \ref{k-keeping}, 
$\langle V(G_1) \cup (S \setminus V(T')) \rangle_{G-E(G_1) \cap E(T')}$ is $k$-connected.
Note that $E(\langle V(G_1) \cup (S \setminus V(T')) \rangle_G) \cap E(T')
= E(G_1) \cap E(T')$.
Thus, we have $\langle V(G_1) \cup (S \setminus V(T')) \rangle_{G-E(G_1) \cap E(T')}
= \langle V(G_1) \cup (S \setminus V(T')) \rangle_{G-E(T')}$. 
Therefore, $\langle V(G_1) \cup (S \setminus V(T')) \rangle_{G-E(T')}$ is $k$-connected. 
For any $w \in V(T') \cap S$, it holds that $|N_{G-E(T')}(w) \setminus V(T')| \geq k$.
Hence, $G-E(T')$ is $k$-connected.
\end{proof}
\bigskip

While fixing the lower bound on $\delta(G)$ in Theorem \ref{Th3}, 
we can show the existence of two disjoint connectivity preserving trees 
when $k = \kappa(G)$. 

\begin{theorem} \label{Th4}
For any trees $T_1$ and $T_2$ of order $m$, every connected cograph $G$
with $\delta(G) \geq \kappa(G)+m-1$ 
contains subtrees $T'_1 \cong T_1$ and $T'_2 \cong T_2$
with $V(T'_1) \cap V'(T'_2) = \emptyset$ 
such that $\kappa(G-E(T'_1) \cup E(T'_2)) = \kappa(G)$. 
\end{theorem}

\begin{proof}
Let $G$ be a connected cograph with $\delta(G) \geq \kappa(G)+m-1$. 
If $m = 1$, then $E(T) = \emptyset$
and the statement vacuously holds. 
Suppose that $m \geq 2$.
Since $\delta(G) \geq \kappa(G)+1$, $G \not\cong K_{n_G}$ 
and the primary cocomponent $G_1$ is disconnected, i.e.,
$G_1$ has at least two components.
Let $G'_1$ and $G''_1$ be components of $G_1$.
Since $\delta(G_1) \geq m-1$, 
$\delta(G'_1) \geq m-1$ and $\delta(G''_1) \geq m-1$.
Thus, $T_1 \cong T'_1 \subseteq G'_1$ and $T_2 \cong T'_2 \subseteq G''_1$
such that $V(T'_1) \cap V'(T'_2) = \emptyset$.
Since $E(T'_1) \cup E(T'_2) \subseteq E(G_1)$, 
from Lemma \ref{edge-deletion},
$\kappa(G-E(T'_1) \cup E(T'_2)) = \kappa(G)$. 
\end{proof}

\section{Variants of Mader's Conjecture for $k$-Edge-Connected Cographs}

In this section, we first show that any connected cograph is maximally edge-connected
and characterize the super edge-connected cographs.

\begin{theorem} \label{super}
Every connected cograph is maximally edge-connected.
Moreover, a connected cograph $G$ is suer edge-connected
if and only if the following two conditions hold:
\begin{enumerate}
\item $G \not\cong C_4$,
\item $G \not\cong  H + K_1$ where $H$ is a disconnected cograph of order at least 4
with $\delta(H) \geq 1$
such that $H$ has a component $H' \cong K_{\delta(G)}$.
\end{enumerate}
\end{theorem}

\begin{proof}
Fig. 4 illustrates the connected cographs of order at most 4.
It can be checked that all of them are maximally edge-connected.
Moreover, $C_4$ is not super edge-connected and 
except for $C_4$, every connected cograph of order at most 4
is super edge-connected.

\begin{figure}[h]
\centering
\epsfig{file=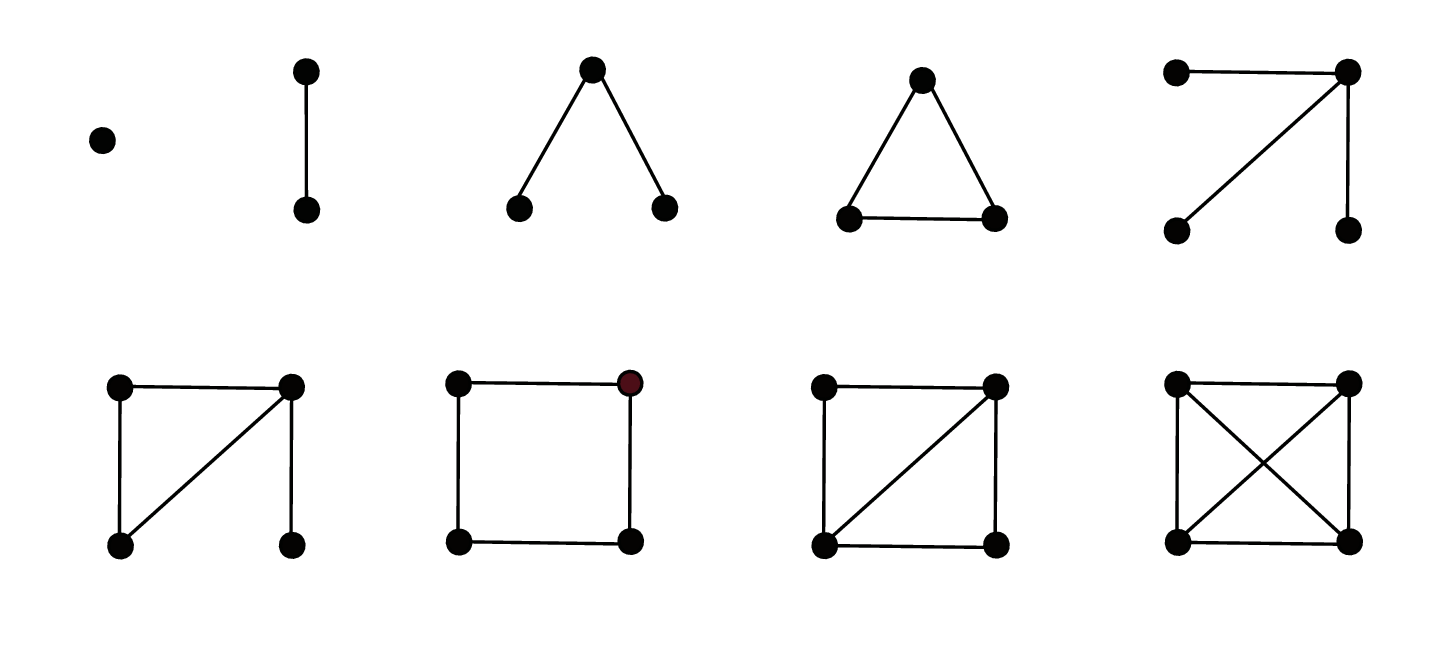,height=45mm}
\caption{The connected cographs of order at most 4.}
\end{figure}

If $G =  G_1 + G_2$ where $G_2 \cong K_1$, $n_{G_1} \geq 4$ and $\delta(G_1) \geq 1$
such that $G_1$ has a component $G'_1 \cong K_{\delta(G)}$,
then the set of $\delta(G)$ edges joining a vertex of $G'_1$ and the vertex of $G_2$
is a minimum edge-cut of $G$ which does not isolate a vertex.
Thus, in such a case, $G$ is not super edge-connected.

Let $G$ be a connected cograph with $n_G \geq 5$. 
Let $S \subsetneq V(G)$ with $2 \leq |S| \leq \lfloor \frac{n_G}{2} \rfloor$.
Note that $|S| \leq |V(G) \setminus S|$.
Let $E_S$ denote the set of edges joining a vertex in $S$ and a vertex in $V(G) \setminus S$,
i.e., $$E_S = \{ uv \in E(G)\ |\ u \in S, v \in V(G) \setminus S \}.$$ 
In order to show that $G$ is maximally edge-connected (respectively, super edge-connected),
it is sufficient to show that $|E_S| \geq \delta(G)$ (respectively, $|E_S| > \delta(G)$). 

Let $S_i = S \cap V(G_i)$ for each $1 \leq i \leq t_G$.
Any vertex $v \in S_i$ is adjacent to every vertex in $\cup_{j \neq i}(V(G_j) \setminus S_j)$
and $|N_{G_i}(v) \setminus S_i|$ vertices in $V(G_i) \setminus S_i$.
Thus, we have 
\begin{equation} \label{Es}
|E_S| = \displaystyle \sum_{i = 1}^{t_G} |S_i|(|V(G) \setminus S|-|V(G_i) \setminus S_i|) + 
\sum_{i = 1}^{t_G} \sum_{v \in S_i}|N_{G_i}(v) \setminus S_i|.
\end{equation}
Let $x \in S_p \subseteq S$ where $1 \leq p \leq t_G$ such that 
${\rm deg}_G(x) = \min_{v \in S}{\rm deg}_G(v)$.
Since
$$\begin{array}{ll}
{\rm deg}_G(x) & = |V(G) \setminus V(G_p)|+{\rm deg}_{G_p}(x) \\
& = |V(G) \setminus V(G_p)|+|N_{G_p}(x) \setminus S_p|+|N_{G_p}(x) \cap S_p| \\
& \leq |V(G) \setminus V(G_p)|+|N_{G_p}(x) \setminus S_p|+|S_p|-1, \\
\end{array}
$$
it holds that 
\begin{equation} \label{Es-degx}
\begin{array}{lll}
|E_S|-{\rm deg}_G(x) 
& \geq & \displaystyle \sum_{i \neq p} |S_i|(|V(G) \setminus S|-|V(G_i) \setminus S_i|) + 
\sum_{i = 1}^{t_G} \sum_{v \in S_i \setminus \{x\}}|N_{G_i}(v) \setminus S_i| \\[5mm]
& & 
+ (|S_p|-1)(|V(G) \setminus S|-|V(G_p) \setminus S_p|)-|S|+1. \\
\end{array}
\end{equation}
Note that $|V(G) \setminus S| - |V(G_p) \setminus S_p| - |V(G) \setminus V(G_p)| - |S_p|
= -|S|$.

We consider two possible cases.

\bigskip
\noindent Case 1: There exists $q \in \{1,2,\ldots,t_G\}$ 
such that $|V(G) \setminus S|-|V(G_q) \setminus S_q| = 0$.

In this case, $V(G) \setminus S = V(G_q) \setminus S_q$ and 
for any $i \neq q$, $S_i = V(G_i)$.
Since 
$$|V(G_q)| \geq |V(G_q) \setminus S_p| = |V(G) \setminus S| \geq |S| 
= |V(G)|-|V(G_q) \setminus S_q| \geq |V(G)| - |V(G_q)|,$$ 
$G_q$ is a primary cocomponent of $G$.
Thus, without loss of generality, we may assume that $q = 1$.
If $p = 1$, then from (\ref{Es-degx}), we have 
$$\begin{array}{lll}
|E_S|-{\rm deg}_G(x) 
& \geq & \displaystyle \sum_{i \neq 1} |S_i| \cdot |V(G) \setminus S| + 
\sum_{v \in S_1 \setminus \{x\}}|N_{G_1}(v) \setminus S_1| -|S|+1.\\[5mm]
& \geq & \displaystyle |S| \left( \sum_{i \neq 1} |V(G_i)| -1 \right)+ 
\sum_{v \in S_1 \setminus \{x\}}|N_{G_1}(v) \setminus S_1|+1 > 0.
\end{array}
$$
Suppose that $p \geq 2$.
Then, from (\ref{Es-degx}), we also have  
$$\begin{array}{lll}
|E_S|-{\rm deg}_G(x) 
& \geq & \displaystyle \sum_{i \neq 1,p} |S_i| \cdot |V(G) \setminus S| + 
\sum_{v \in S_1}|N_{G_1}(v) \setminus S_1| \\[5mm]
& & 
+ (|S_p|-1)|V(G) \setminus S|-|S|+1.\\[5mm]
& \geq & \displaystyle |S| \left( \sum_{i \neq 1} |V(G_i)| -2 \right)
+\sum_{v \in S_1}|N_{G_1}(v) \setminus S_1|  + 1. 
\end{array}
$$
Thus, if $t_G \geq 3$, or $t_G = 2$ and $|V(G_2)| \geq 2$, then $|E_S|-{\rm deg}_G(x) > 0$. 
Assume that $t_G = 2$ and $|V(G_2)| = 1$, i.e., $V(G_2) = \{x\}$.
Then, ${\rm deg}_G(x) = n_G-1$ and for any $v \in S$, ${\rm deg}_G(v) = n_G-1$.
Since $S_1 \neq \emptyset$, $G_1$ must be a nontrivial connected cograph, 
which is a contradiction. 

Hence, in Case 1, it holds that $|E_S| > {\rm deg}_G(x) \geq \delta(G)$. 

\bigskip
\noindent Case 2: For all $1 \leq i \leq t_G$, $|V(G) \setminus S|-|V(G_i) \setminus S_i| \geq 1$. 

In this case, from (\ref{Es-degx}), we have 
$$|E_S|-{\rm deg}_{G}(x) \geq 
\displaystyle \sum_{i = 1}^{t_G} \sum_{v \in S_i \setminus \{x\}}|N_{G_i}(v) \setminus S_i| \geq 0.$$ 
This means that $|E_S| \geq \delta(G)$.
At this juncture, our proof for maximally edge-connectedness is complete. 

Two cases are considered, depending on $t_G$. 

\bigskip
\noindent Case 2.1: $t_G \geq 3$.

Since $|V(G) \setminus S|-|V(G_1) \setminus S_1| \geq 1$,
there exists an $i_1 \in \{2,\ldots,t_G\}$ such that
$V(G_{i_1}) \setminus S_{i_1} \neq \emptyset$.
Then, there exists an $i_2 \in \{1,2,\ldots,t_G\} \setminus \{i_1\}$ such that
$V(G_{i_2}) \setminus S_{i_2} \neq \emptyset$. 
Note that for any $i \in \{1,2,\ldots,t_G\} \setminus \{i_1,i_2\}$,
$|V(G) \setminus S|-|V(G_i) \setminus S_i| \geq 2$.
If there exists a $j \in \{1,2,\ldots,t_G\} \setminus \{i_1,i_2\}$ such that 
$S_j \setminus \{x\} \neq \emptyset$, then from (\ref{Es-degx}), 
we have $|E_S|-{\rm deg}_G(x) \geq |S_j \setminus \{x\}| > 0$.
If there exists a $j \in \{1,2,\ldots,t_G\} \setminus \{i_1,i_2\}$ such that 
$V(G_j) \setminus S_j \neq \emptyset$, then
for all $1 \leq i \leq t_G$, 
$|V(G) \setminus S|-|V(G_i) \setminus S_i| \geq 2$
and from (\ref{Es-degx}), we have $|E_S|-{\rm deg}_G(x) \geq |S \setminus \{x\}| > 0$.
Assume that for any $j \in \{1,2,\ldots,t_G\} \setminus \{i_1,i_2\}$,
$S_j \setminus \{x\} = \emptyset$ and $V(G_j) \setminus S_j = \emptyset$.
Then, it follows that $t_G = 3$, $\{i_1,i_2,p\} = \{1,2,3\}$
and $V(G_p) = \{x\}$. 
Thus, we have ${\rm deg}_G(x) = n_G-1$ and 
for any $v \in S_{i_1} \cup S_{i_2}$, ${\rm deg}_G(v) = n_G-1$.
Note that $S_{i_1} \cup S_{i_2} \neq \emptyset$, 
$V(G_{i_1}) \setminus S_{i_1} \neq \emptyset$ and 
$V(G_{i_2}) \setminus S_{i_1} \neq \emptyset$.
Therefore, $G_{i_1}$ or $G_{i_2}$ must be a nontrivial connected cograph, which is contradictory.

\bigskip
\noindent Case 2.2: $t_G = 2$.

Since $|V(G) \setminus S| \geq \lceil \frac{n_G}{2} \rceil \geq 3$, 
$|V(G_1) \setminus S_1| \geq 2$ or 
$|V(G_2) \setminus S_2| \geq 2$.
Note that $V(G_1) \setminus S_1 \neq \emptyset$ and
$V(G_2) \setminus S_2 \neq \emptyset$. 
If $|V(G_1) \setminus S_1| \geq 2$ and $|V(G_2) \setminus S_2| \geq 2$, then we have 
$|E_S|-{\rm deg}_G(x) > 0$.
If $|V(G_1) \setminus S_1| = 1$ and $|V(G_2) \setminus S_2| \geq 2$, 
then $|S_1| \geq 2$, i.e., $|S_1 \setminus \{x\}| \geq 1$ 
and thus we also have $|E_S|-{\rm deg}_G(x) > 0$.

Suppose that $|V(G_1) \setminus S_1| \geq 2$ and $|V(G_2) \setminus S_2| = 1$.
If $|S_2 \setminus \{x\}| \geq 1$, then $|E_S|-{\rm deg}_G(x) > 0$.
We then consider two remaining cases.

\bigskip
\noindent Case 2.2.1: $S_2 = \{x\}$.

In this case, $|V(G_2)| = 2$, i.e., $G_2 \cong \overline{K_2}$.
From (\ref{Es}), we have $|E_S| \geq |S_1|+|V(G_1) \setminus S_1| = |V(G_1)| = n_G-2 \geq 3$.
Since ${\rm deg}_G(x) = n_G-2$, for any $v \in S_1$, ${\rm deg}_G(v) \geq n_G-2$, i.e., 
${\rm deg}_{G_1}(v) \geq |V(G_1)|-2$. 
Note that $S_1 \neq \emptyset$.
Since $G_1$ is disconnected, $G_1$ has an isolated vertex,
which implies that $\delta(G) = 2$.
Therefore, $|E_S|  > \delta(G)$.

\bigskip
\noindent Case 2.2.2: $S_2 = \emptyset$.

In this case, $G_2 \cong K_1$.
From (\ref{Es}), we have
$|E_S| = |S_1|+ \sum_{v \in S_1}|N_{G_1}(v) \setminus S_1|$.
Since ${\rm deg}_G(x) = |N_{G_1}(x) \setminus S_1|+|N_{G_1}(x) \cap S_1|+1$, 
if $|N_{G_1}(x) \cap S_1| < |S_1|-1$, then
$|E_S| > {\rm deg}_G(x)$.
Moreover, if there exists a vertex $y \in S_1 \setminus \{x\}$ such that
$N_{G_1}(y) \setminus S_1 \neq \emptyset$, then 
$|E_S| > {\rm deg}_G(x)$.
From these observations and the fact that ${\rm deg}_G(x) = \min_{v \in S}{\rm deg}_G(x)$,
it follows that $|E_S| > \delta(G)$, except for the case that 
$\langle S_1 \rangle_{G_1}$ is a component of $G_1$ such that 
$\langle S_1 \rangle_{G_1} \cong K_{|S_1|}$ where $|S_1| = {\rm deg}_G(x) = \delta(G)$. 
Note that when ${\rm deg}_G(x) > \delta(G)$, we have the desired result
even if $|E_S| = {\rm deg}_G(x)$. 

Therefore, in Case 2, it holds that $|E_S| \geq \delta(G)$, and
$|E_S| > \delta(G)$ if $G$ satisfies the second condition in the statement. 
\end{proof}
\bigskip

From Theorem \ref{super}, we can see that 
if a connected cograph $G \not\cong C_4$ is not super edge-connected, then
$\kappa(G) = 1$.
Thus, we have the following. 

\begin{corollary} \label{super-cor}
Every 2-connected cograph $G \not\cong C_4$ is super edge-connected.
\end{corollary}

While any maximally connected graph is maximally edge-connected and so is 
any super edge-connected graph, 
there is no inclusion relationship between the maximally connected garphs and
the super edge-connected graphs in general. 
However, from Theorems \ref{maximal} and \ref{super}, the following result is obtained. 

\begin{corollary}
Every maximally connected cograph $G \not\cong C_4$ is super edge-connected.
\end{corollary}

From Corollary \ref{Co1} for $k = 1$ and Theorem \ref{super}, 
we have the following statement as a corollary. 

\begin{corollary} \label{maximal-edge-keeping}
For any tree $T$ of order $m$, every maximally edge-connected cograph $G$
with $\delta(G) \geq m$ contains a subtree $T' \cong T$ 
such that $G-V(T')$ is maximally edge-connected. 
\end{corollary}

Define the cograph $G$ as follows:
$$G = (K_{m-1} \cup K_{m-1} \cup K_{m-1}) + K_1.$$
Then, $G$ is a maximally edge-connected cograph with $\delta(G) = m-1$.
For any subtree $T'$ of order $m$ in $G$,
$G-V(T')$ is disconnected.
Thus, the lower bound of $m$ on $\delta(G)$ in Corollary \ref{maximal-edge-keeping}
is tight in the sense that $G-V(T')$ is connected. 

Applying Theorem \ref{super},
we show the following result on an edge-connectivity keeping tree, 
where the lower bound of $k+m-[k = 1]$ on $\delta(G)$ is tight, 
which follows from the case that $G \cong K_{k+m-[k = 1]}$. 

\begin{theorem} \label{Th6}
For any tree $T$ of order $m$, every $k$-edge-connected cograph $G$
with $\delta(G) \geq k+m-[k = 1]$ 
contains a subtree $T' \cong T$ such that $G-V(T')$ is $k$-edge-connected.
\end{theorem}

\begin{proof}
Let $T$ be a tree of order $m$
and $G$ a $k$-edge-connected cograph with $\delta(G) \geq k+m-[k = 1]$. 
When $k = 1$, the statement is essentially the same as Corollary \ref{Co1}.
Note that any nontrivial graph is 1-edge-connected if and only if it is 1-connected,
and the trivial graph is assumed to be 1-connected and 1-edge-connected. 
Suppose that $k \geq 2$.
Then, from Corollary \ref{Co1} for $k = 1$, 
$G$ contains $T' \cong T$ such that $G-V(T')$ is connected.
Since $\delta(G-V(T')) \geq k$ and by Theorem \ref{super}, 
any connected cograph is maximally edge-connected, 
$G-V(T')$ is $k$-edge-connected. 
Hence, the statement also holds for $k \geq 2$. 
\end{proof}
\bigskip

We next present a result on a super edge-connectedness keeping tree. 

\begin{theorem} \label{edge-keeping}
For any tree $T$ of order $m$, every super edge-connected cograph $G$
with $\delta(G) \geq m+2$
contains a subtree $T' \cong T$ such that $G-V(T')$ is super edge-connected.
\end{theorem}

\begin{proof}
Let $T$ be a tree of order $m$
and $G$ a super edge-connected cograph with $\delta(G) \geq m+2$.
From Theorem \ref{super}, $G \not\cong C_4$ and 
$G \not\cong H + K_1$ where $H$ is a disconnected cograph with $n_{H} \geq 4$ 
and $\delta(H) \geq 1$ 
such that $H$ has a component $H' \cong K_{\delta(G)}$.

If $|V(G_1)| = 1$, then $G \cong K_{n_G}$ and for any $T \cong T' \subset G$,
$G-V(T') \cong K_{n_G-m}$; thus $G-V(T')$ is super edge-connected. 
Suppose that $|V(G_1)| \geq 2$ and $n_G \geq m+4$. 
Let $S = V(G) \setminus V(G_1)$. 

We consider two cases depending on $|S|$. 

\bigskip
\noindent Case 1: $|S| \geq 2$.

Since $G$ is 2-connected and $\delta(G) \geq m+2$, from Corollary \ref{Co1} for $k = 2$,
there exists a subtree $T' \cong T$ of $G$ such that $G-V(T')$ is 2-connected; 
in more detail, from Cases 2, 3, and 4 in the proof of Theorem \ref{Th1},
$|V(G_1) \setminus V(T')| \geq 2$ and $|S \setminus V(T')| \geq 2$.
Thus, by Corollary \ref{super-cor}, if $G-V(T') \not\cong C_4$, then 
$G-V(T')$ is super edge-connected.

Suppose that $G-V(T') \cong C_4$, i.e., 
$\langle V(G_1) \setminus V(T') \rangle_G \cong \langle S \setminus V(T') \rangle_G \cong
\overline{K_2}$.
Since $\delta(G) \geq m+2$, every vertex of $G-V(T')$ is adjacent to every vertex of $T'$.
Thus, if $V(G_1) \cap V(T') \neq \emptyset$, then
$G_1$ is nontrivial and connected, which is a contradiction.
Therefore, $V(G_1) \cap V(T') = \emptyset$, i.e., $V(T') \subset S$.
Let $z$ be a leaf of $T'$ such that $wz \in E(T')$.
Let $T''$ be a subtree obtained from $T'-z$ by adding a vertex $x \in V(G_1)$
with the edge $wx$. 
Then, $T'' \cong T'$ such that $|V(G_1) \setminus V(T'')| = 1$
and $|S \setminus V(T'')| = 3$, i.e., $G-V(T'') \not\cong C_4$.
Therefore, $G-V(T'')$ is super edge-connected.

\bigskip
\noindent Case 2: $|S| = 1$.

In this case, $G = G_1 + K_1$.
If every component of $G_1$ is a complete graph, then a component of the least order 
must be isomorphic to $K_{\delta(G)}$, which contradicts the assumption of $G$. 
Thus, $G_1$ has a component $G'_1$ which is not a complete graph.
Note that $G'_1$ is a connected cograph with $\delta(G'_1) \geq \delta(G)-1 \geq m+1$. 
Let $G'_{1,1}$ be a primary cocomponent of $G'_1$.
Then, $|V(G'_{1,1})| \geq 2$ and $G'_{1,1}$ is disconnected.
Let $x', y' \in V(G'_{1,1})$ such that $x'y' \not\in E(G)$.
Also let $S' = V(G'_1) \setminus V(G'_{1,1})$.

Now assume that there is a subtree $T^\ast \cong T$ in $G'_1$
such that $G'_1-V(T^\ast)$ is a connected graph but not a complete graph.
Then, $G_1 - V(T^\ast)$ is disconnected and
$(G_1 - V(T^\ast))+K_1 = G-V(T^\ast)$, i.e., $\kappa(G-V(T^\ast)) = 1$.  
Thus, we have $G-V(T^\ast) \not\cong C_4$. 
Since $\delta(G) \geq \delta(G-V(T^\ast))$ and
$G_1$ has no component isomorphic to $K_{\delta(G)}$, 
$G_1-V(T^\ast)$ also has no component isomorphic to $K_{\delta(G-V(T^\ast))}$.
Thus, $G-V(T^\ast)$ is super edge-connected.
Therefore, it is sufficient to show the existence of such a subtree $T^\ast$ in $G'_1$.
We distinguish two cases. 

\bigskip
\noindent Case 2.1: $E(\langle S' \rangle_G) \neq \emptyset$.

Let $w'z' \in E(\langle S' \rangle_G)$.
Since $\delta(G'_1-\{x',y'\}) \geq m-1$,
by Lemma \ref{tree2},
we can construct a subtree $T' \cong T$ of $G'_1 - \{x',y'\}$ 
so that $w'z' \in E(T')$ and $z'$ is a leaf.
If $S' \setminus V(T') \neq \emptyset$, then 
$G'_1-V(T')$ is a connected graph but not a complete graph. 
Suppose that $S' \setminus V(T') = \emptyset$.
Since $\delta(G'_1) \geq m+1$ and $G'_1$ is not a complete graph, 
$|V(G'_{1})| \geq m+3$.
Thus, $|V(G'_{1,1}) \setminus V(T')| \geq 3$. 
Let $T''$ be a subtree obtained from $T'-z'$ by adding a vertex 
$z'' \in V(G'_{1,1}) \setminus (V(T') \cup \{x',y'\})$ with the edge $w'z''$.
Then, $T''$ is a desired subtree of $G'_1$.

\bigskip
\noindent Case 2.2: $E(\langle S' \rangle_G) = \emptyset$. 

Suppose that $|S'| \geq 2$.
If $E(G'_{1,1}) \neq \emptyset$, then
we can apply a similar discussion for Case 2.1 
by selecting two vertices in $S'$ and an edge of $G'_{1,1}$ 
instead of $x',y' \in V(G'_{1,1})$ and $w'z' \in E(\langle S' \rangle_G)$.
Thus, in such a case, we have the desired result.
Suppose that $E(G'_{1,1}) = \emptyset$.
Then $G'_{1} \cong K_{|V(G'_{1,1})|,|S'|}$.
Since $\delta(G'_1) \geq m+1$, $|V(G'_{1,1})| \geq m+1$ and $|S'| \geq m+1$.
Thus, $G'_1$ contains a subtree $T' \cong T$ such that 
$G'_1 - V(T') \cong K_{n_1,n_2}$, where $n_1 \geq 2$ and $n_2 \geq 2$.
Therefore, $T'$ is a desired subtree of $G'_1$.

Suppose that $|S'| = 1$.
Let $G''_{1,1}$ be a component of $G'_{1,1}$.
Since $\delta(G''_{1,1}) \geq m$, there exists a subtree $T'' \cong T$ in $G''_{1,1}$
such that $|V(G''_{1,1}-V(T''))| \geq 1$. 
Therefore, $G'_{1,1}-V(T'')$ is disconnected
and $(G'_{1,1}-V(T''))+K_1 = G'_1 - V(T'')$ is a connected graph but not a complete graph.
Hence, $T''$ is a desired subtree of $G'_1$. 
\end{proof}

\bigskip

Let $m$ be even. 
Define the cograph $G$ as $G = (H_1 \cup H_2) + K_1$, 
where 
$$H_1 \cong H_2 \cong 
\underbrace{(K_1 \cup K_1) + (K_1 \cup K_1) + \cdots + (K_1 \cup K_1)}_{\frac{m}{2}+1}.$$
Then, $G$ is a super edge-connected cograph with $\delta(G) = m+1$. 
For any subtree $T' \subset G$ isomorphic to the star of order $m$, 
$G-V(T')$ is a disconnected graph with $\delta(G-V(T')) \geq 1$, 
or $G-V(T') \cong (K_2 \cup H_2) + K_1$; thus, in either case, 
$G-V(T')$ is not super edge-connected.
Hence,  the lower bound of $m+2$ on $\delta(G)$ in Theorem \ref{edge-keeping} is tight.

Given a nontrivial tree $T$ with bipartition $(V_1,V_2)$, 
let $\beta(T) = \min\{|V_1|, |V_2|\}$.
We then show the following result on an edge-connectivity preserving tree. 

\begin{theorem} \label{Th5}
For any tree $T$ of order $m$, every $k$-edge-connected cograph $G$
with $\delta(G) \geq \max\{k+\Delta(T)+\beta(T)-1, m-1\}$ 
contains a subtree $T' \cong T$ such that $G-E(T')$ is $k$-edge-connected. 
\end{theorem}

\begin{proof}
When $m = 1$, the statement vacuously holds since $E(T) = \emptyset$. 
Let $T$ be a nontrivial tree of order $m$ with bipartition $(V_1,V_2)$, where $|V_1| \geq |V_2|$.
Let $G$ be a $k$-edge-connected cograph of order $n$ with 
$\delta(G) \geq \max\{k+\Delta(T)+|V_2|-1, m-1\}$.
Since $\delta(G) \geq m-1$, $G$ contains a subtree $T' \cong T$ with 
bipartition $(V'_1,V'_2)$, where $|V'_1| \geq |V'_2|$.
Note that $E(G-V'_2) \cap E(T') = \emptyset$. 

We employ induction on the order $n$ of $G$.
Consider the base case that $n = \max\{k+\Delta(T)+|V_2|, m\}$, i.e., $G \cong K_{n}$.
Since $G-V'_2 \cong K_{n-|V'_2|}$, $\kappa(G-V'_2) \geq k+\Delta(T)-1 \geq k$.
Any $v \in V'_2$ is adjacent to at least  $|V(G) \setminus V'_2|-\Delta(T)$ vertices in
$V(G) \setminus V'_2$ through edges in $E(G) \setminus E(T')$, where 
$|V(G) \setminus V'_2|-\Delta(T) \geq k$.
Therefore, 
$G-E(T')$ is $k$-connected; thus it is $k$-edge-connected.

Suppose that $n >  \max\{k+\Delta(T)+|V_2|, m\}$.
We consider two cases depending on whether $G-V'_2$ is connected or not. 

\bigskip
\noindent Case 1: $G-V'_2$ is connected.

By Theorem \ref{super}, $G-V'_2$ is maximally edge-connected.
Since $\delta(G-V'_2) \geq k+\Delta(T)-1 \geq k$,
$G-V'_2$ is $k$-edge-connected.
For any $w \in V'_2$, 
$$\begin{array}{ll}
|N_{G-E(T')}(v) \cap (V(G) \setminus V'_2)| 
& = |N_{G-E(T')}(v) \setminus V'_2| \\
& \geq {\rm deg}_G(v) - {\rm deg}_{T'}(v) - (|V'_2|-1)  \\
& \geq \delta(G)-\Delta(T)-|V_2|+1 \\ 
& \geq k.
\end{array}$$
Thus, $G-E(T')$ is $k$-edge-connected. 

\bigskip
\noindent Case 2: $G-V'_2$ is disconnected.

In this case, there exists a nontrivial cocomponent $G_j$ 
such that $V(G) \setminus V(G_j) \subseteq V'_2$.  
Thus, $V'_1 \subseteq V(G_j)$. 
Since 
$$|V(G_j)| \geq |V'_1| \geq |V'_2| \geq |V(G) \setminus V(G_j)|,$$ 
$G_j$ is a primary cocomponent of $G$.
Without loss of generality, 
we may assume that $j = 1$.
Let $S = V(G) \setminus V(G_1)$.
Then, $V'_1 \subseteq V(G_1)$, $S \subseteq V'_2$ and $|V(G_1)| \geq |S|$.
We distinguish two cases. 

\bigskip 
\noindent Case 2.1: $E(\langle S \rangle_G) \neq \emptyset$.

By Lemma \ref{tree2},
we can construct $T'' \cong T$ from an edge $uv \in E(\langle S \rangle_G)$.
Let $(V''_1,V''_2)$ be the bipartition of $T''$ where $|V''_1| \geq |V''_2|$.
Without loss of generality, we may suppose that $u \in V''_1$ and $v \in V''_2$.
Now assume that $G-V''_2$ is disconnected.
Since $u \in V''_1 \setminus V(G_1)$, 
$V''_1$ must be included in another primary cocomponent of $G_\ell$ where $\ell \neq 1$,
i.e., $V''_1 \subseteq V(G_\ell)$ and $V(G) \setminus V(G_\ell) \subseteq V''_2$.
Since $V(G_1) \subseteq V(G) \setminus V(G_\ell)$ and $v \in V''_2 \setminus V(G_1)$, 
we have $$|V''_2| > |V(G_1)| \geq |S| \geq |V(G_\ell)| \geq |V''_1|,$$ which is a contradiction. 
Thus, $G-V''_2$ is connected.
Therefore, similarly to Case 1, we obtain the desired result.

\bigskip
\noindent Case 2.2: $E(\langle S \rangle_G) = \emptyset$.

In this case, $\langle S \rangle_G$ is a cocomponent of $G$, i.e., 
$\langle S \rangle_G = G_2 \cong \overline{K_{|S|}}$ and $t_G = 2$. 
Note that for any $v \in S$, it holds that ${\rm deg}_G(v) = |V(G_1)| \geq \delta(G) \geq m-1$. 
We finally consider two cases depending on $|S|$. 

\bigskip
\noindent Case 2.2.1: $|S| \geq 2$.

Let $x$ be a leaf of $T$ such that $xy \in E(T)$.
Let $w,z \in S$.
Since $\delta(G-w) \geq m-2$, by Lemma \ref{tree2}, 
we can construct $T'_x \cong T-x$ in $G-w$
so that $z \in V(T'_x)$ and $z$ corresponds to $y$ of $T-x$.
Since $z$ is adjacent to every vertex in $V(G_1)$, 
we can extend $T'_x$ to a subtree $T^\ast \cong T$ in $G-w$
by adding a vertex $z' \in V(G_1) \setminus V(T'_x)$ with the edge $zz'$. 
Let $(V^\ast_1,V^\ast_2)$ be the bipartition of $T^\ast$ where $|V^\ast_1| \geq |V^\ast_2|$.
If $V(G_1) \setminus V^\ast_2 \neq \emptyset$, then
$G-V^\ast_2$ is connected since $w \not\in V(T^\ast)$. 
Thus, in such case, we have the desired result similarly to Case 1.
If $V(G_1) \subseteq V^\ast_2$, then $m-1 \leq |V(G_1)| \leq |V^\ast_2| \leq \frac{m}{2}$, i.e., 
$m = 2$ and $|V(G_1)| = 1$, which contradicts the facts that $|V(G_1)| \geq |S| \geq 2$. 

\bigskip
\noindent Case 2.2.2: $|S| = 1$.

Let $S = \{w\}$.
Let $G'_1$ be a component of $G_1$ and $W_1 = \langle V(G'_1) \cup \{w\} \rangle_G$.
Note that $\delta(W_1) \geq \delta(G)$ and $W_1$ is a $k$-edge-connected cograph.
By the inductive hypothesis, 
$W_1$ contains a subtree $T'' \cong T$ such that $W_1 - E(T'')$ is $k$-edge-connected. 
For any other component $G'_i$ of $G_1$ where $i \neq 1$,
$\langle V(G'_i) \cup \{w\} \rangle_G$ is also $k$-edge-connected.
Thus, $G-E(T'')$ is the graph obtained from disjoint $k$-edge-connected graphs by 
identifying one vertex.
Hence, $G-E(T'')$ is $k$-edge-connected.
\end{proof}
\bigskip

Since $k+\Delta(T)+\beta(T)-1 \leq k+m-1$, 
Conjecture \ref{H1} holds for $k$-edge-connected cographs.
The tightness of the lower bound of $k+m-1$ on $\delta(G)$ 
can be checked similarly to Theorem \ref{Th3}. 

\begin{corollary} 
For any tree $T$ of order $m$, every $k$-edge-connected cograph $G$
with $\delta(G) \geq k+m-1$ 
contains a subtree $T' \cong T$ such that $G-E(T')$ is $k$-edge-connected.
\end{corollary}

\section{Concluding Remarks}

In this paper, we have shown that Mader's conjecture is true for $k$-connected cographs.
We also show that three variants of Mader's conjecture hold for $k$-connected or
$k$-edge-connected cographs.
Our proofs are all constructive and lead to algorithms for finding desired trees. 
While the lower bounds on the minimum degree in the variants of Mader's conjecture
are tight even for cographs, 
it remains unknown whether the lower bound on the minimum degree in Mader's conjecture 
can be improved or not for cographs. 

We have furthermore presented tight lower bounds on the minimum degree of a cograph for
the existence of two disjoint connectivity keeping trees, 
a maximal connectedness keeping tree 
and a super edge-connectedness keeping tree,
and have mentioned a result on an ideal connectedness keeping tree.
It would be interesting to investigate these variants of a connectivity keeping tree 
for other graph classes or general graphs. 

As far as we know, the class of cographs is the first graph class for which Mader's conjecture 
and its variants hold for all $k \geq 1$, except for specific graphs such as complete graphs
and paths.
Apart from Mader's conjecture, the Fox-Sudakov conjecture have recently been shown to be true
for cographs \cite{FNSS}.
The class of cographs is equivalent to the class of $P_4$-free graphs.
It would also be interesting to study Mader's conjecture and its variants for a larger graph class
such as the class of $P_\ell$-free graphs where $\ell \geq 5$. 
In fact, another conjecture, the Erd\H{o}s-Hajnal conjecture
have been investigated for $P_5$-free graphs \cite{BB}.

\end{document}